\documentclass[reqno]{amsart}
\usepackage{amsfonts,amsmath,latexsym,verbatim,amscd,mathrsfs,xcolor,array, upgreek}
\usepackage{amssymb,amsthm,graphicx,color}
\usepackage[hmargin=2.5cm, vmargin=2.5cm]{geometry}
\usepackage{epstopdf}
\usepackage{cite}
\usepackage[english]{babel}
\usepackage{multirow}
\usepackage[font=bf,aboveskip=15pt]{caption}
\usepackage[toc,page]{appendix}
\usepackage[colorlinks, linkcolor=black, citecolor=blue]{hyperref}

\usepackage{enumerate}
\usepackage{bookmark}
\usepackage{float}
\usepackage{mathtools}
\usepackage{bm}

\newcommand{\dd}{ \mathrm{d}}

\newcommand{\R} {\mathbb R}

\newcommand{\cuad}{{\sqcap\kern-.68em\sqcup}}

\newcommand{\dist}{{\rm dist}  }

\newcommand{\be}{\begin{equation}}
	\newcommand{\ee}{\end{equation}}

\newcommand{\bv}[1]{\boldsymbol{#1}}  
\newcommand{\eps}{\varepsilon}
\newcommand{\M}{\mathcal M}
\newcommand{\Dif}[1]{\mathcal D[#1]}   
\newcommand{\Kop}{\mathcal K}   
\newcommand{\ap}{\mathrm{ap}}
\newcommand{\nY}[1]{\|#1\|_{\mathcal Y}}

\newcommand{\Lin}{\mathcal L_\beta}
\newcommand{\nX}[1]{\|#1\|_{\mathcal X}}
\newcommand{\qbar}{\bar q} 
\newcommand{\D}{\mathcal D_n}
\newcommand{\rad}{\mathrm{rad}}
\newcommand{\bx}[1]{\boldsymbol{#1}}
\newcommand{\sgn}{\operatorname{sgn}}
\newcommand{\Lr}{\mathcal L}

\newtheorem{theorem}{Theorem}[section]
\newtheorem{definition}[theorem]{Definition}
\newtheorem{lemma}[theorem]{Lemma}
\newtheorem{proposition}[theorem]{Proposition}

\newtheorem{corollary}[theorem]{Corollary}
\newtheorem{remark}{Remark}[section]
\newcommand{\bremark}{\begin{remark} \em}
	\newcommand{\eremark}{\end{remark} }

\newcommand{\revised}[1]{{\color{black}#1}}
\numberwithin{equation}{section}

\begin{document}
	\title[Multiple solutions for competition systems]{\revised{Multiple positive solutions near segregated profiles in Lotka-Volterra competition systems with asymmetric interactions}}
	\author[Z. Li]{Zaizheng Li}
		\address{\noindent
		School of Mathematical Sciences, Hebei Normal University, Shijiazhuang 050024, Hebei, P. R. China}
	\email{zaizhengli@hebtu.edu.cn}
	
	\author[S. Terracini]{Susanna Terracini}
	\address{\noindent
		Dipartimento di Matematica “Giuseppe Peano”, Universit\`a di Torino, Via Carlo Alberto, 10, 10123 Torino, Italy}
	\email{susanna.terracini@unito.it}

%	\tableofcontents
	
	\begin{abstract}
\revised{
We consider stationary Lotka--Volterra competition-diffusion systems with any number
of species and positive, possibly asymmetric competition coefficients. We
study two settings with homogeneous Dirichlet boundary conditions: a
one-dimensional interval with constant drift, and radial solutions on a ball.
In each setting we first determine when a prescribed ordering of species admits a
segregated profile. Starting from such a profile, we construct a positive
solution for every sufficiently large competition parameter $\beta$. The
uniform distance from the profile is of order $\varepsilon=\beta^{-1/3}$, and the
rescaled densities converge near each interface to a two-species model.
Counting the admissible orderings gives a lower bound on the number of
positive solutions. We also prove a non-approximation result when one
directed competition coefficient is zero and the other species are
suppressed on the interval or shell used in the proof.}

	\end{abstract}
	
	%\keywords{Competing species, Segregated profile, Converse problem,  Perturbation argument}
	
	%%\pacs[JEL Classification]{D8, H51}
%	\pacs[MSC Classification (2020)]{35K51, 35B32}
%\subjclass[2020]{Primary 35B40, 35B36; Secondary 35B20}
%	\noindent\footnotesize	{ \textbf{AMS Subject Classification (2020):} 35K51, 35B32, 92D25.}
		\maketitle
	
	\section{Introduction}
	
	In this article, we consider pattern formation for the following Gause-Lotka-Volterra competition-diffusion system:
	\begin{equation}\label{Volterra-system}
		\partial_t u_i-\Delta u_i
		=f_i(u_i)-\beta u_i\sum_{j\neq i}\alpha_{ij}u_j,
		\quad i=1,2,\ldots,k.
	\end{equation}
	  Here, $u_i$ denotes the population density of the $i$-th species, $\Delta u_i$ describes its spatial dispersal, and $f_i(u_i)$ represents its intrinsic population growth. The parameter $\beta>0$ measures the overall strength of interspecific competition, while the coefficients $\alpha_{ij}>0$ ($i\neq j$) characterize the competitive interactions between different species and are not necessarily symmetric (see e.g. \cite{Cantrell-Cosner-2003}). 
	 % In this context, $f_i(u)$ are typically of logistic type: $f_i(s)=\mu s(1-s)$.
	 
	  As $\beta\to+\infty$, the strong competition between species leads to spatial segregation, that is,  distinct species cannot coexist at the same point, as proved in \cite{Dancer-Du-1994, 2005CTV, Wang-Zhang-2010, Dancer-Wang-Zhang-2012}. Consequently, the limiting segregated profiles have pairwise disjoint supports, separated by free boundaries. 
	  
	  Related work on two-species systems addresses competitive
exclusion and coexistence
\cite{Matano-Mimura-1983,Cosner-Lazer-1984,Dancer-1984,Dancer-1991,Gui-Lou-1994}
and segregation in the time-dependent problem
\cite{Dancer-Hilhorst-Mimura-Peletier-1999,Dancer-Zhang-2002,Crooks-Dancer-Hilhorst-Mimura-Ninomiya-2004},
while studies of rotating waves address their existence and the
selection of their angular speed
\cite{Salort-Terracini-Verzini-Zilio-2025,Li-Terracini-2025,spadaro2026selection}.

	 In this paper we are interested in the {\it converse problem}: \revised{given a segregated profile, can we construct positive solutions  of the competition system for sufficiently large but finite $\beta$ that converge
to it as $\beta$ tends to infinity? }

%	 Before addressing the converse problem, we first review some classical results %concerning the existence, stability, uniqueness, multiplicity, and the asymptotic behavior of %positive solutions to system \eqref{Volterra-system}.

\revised{We consider stationary solutions, allowing some constant drift. For two species, the competition terms cancel in the equation for
$\alpha_{21}u_1-\alpha_{12}u_2$. Dancer and Du \cite{Dancer-Du-1994}
used this cancellation and topological degree to construct solutions from
a nondegenerate limit profile with logistic growth. Conti, Terracini and
Verzini \cite{2005CTV} obtained quantitative estimates for convergence to
segregated limits in the two-species case. For three or more species, the
existence results in \cite{Dancer-Du-1995,Dancer-Du-1995-2} concern positive
steady states in suitable parameter ranges. Other results concern the
regularity of solutions and segregated limits
\cite{2005CTV,Soave-Zilio-2015,Zhang-2026} and the structure of the
interfaces
\cite{Conti-Terracini-Verzini-2006,Wang-Zhang-2010,Terracini-Verzini-Zilio-2019}.
}

\revised{The converse problem also arises for two-component
Gross--Pitaevskii systems, that possess a variational structure. On bounded domains, the results of Noris,
Tavares, Terracini and Verzini
\cite[Corollary~1.6]{Noris-Tavares-Terracini-Verzini-2012},
combined with those of Dancer, Wang and Zhang
\cite{Dancer-Wang-Zhang-2012}, give solutions converging to
segregated profiles in $H^1\cap C^{0,\alpha}$. Dancer
\cite{Dancer-2014} used a mountain-pass argument to obtain solutions
converging to a least-energy segregated profile for a wider range of
parameters.

The construction of Casteras and Sourdis
\cite{Casteras-Sourdis-2020,casteras2017converse} is closer to the
one we use here. Starting from a nondegenerate radial profile on a ball
or an annulus, they built an approximate solution near the interface
and corrected it to obtain an exact solution with a quantitative
convergence rate. Kowalczyk, Pistoia and Vaira
\cite{Kowalczyk-Pistoia-Vaira-2023} extended this approach to bounded
planar domains under suitable nondegeneracy assumptions. Aftalion and
Sourdis \cite{Aftalion-Sourdis-2017} used a related perturbation argument
on the whole line. In the whole space, Esposito, Figueroa, Pistoia and
Vaira \cite{esposito2026entire} constructed infinitely many nonradial
positive solutions for systems with power-law interactions of the form
$\beta u_1^{a_1}u_2^{a_2}$.

We study the converse problem for Lotka--Volterra systems
with $k\ge2$ species on an interval with constant drift and in the
radial setting on a ball. For each segregated profile described below,
we construct positive solutions for all sufficiently large $\beta$.
Their uniform distance from the prescribed profile is of order
$\beta^{-1/3}$.

Our proof follows the perturbation argument of Casteras and Sourdis
\cite{Casteras-Sourdis-2020}. At each interface, only two species meet,
so the mutual competition terms cancel in a weighted difference of
their equations. To use this cancellation in the linearized problem,
we need nondegeneracy both on each occupied interval or shell and
across the interfaces. Casteras and Sourdis assume these two
properties in \cite[Theorem~1.1]{Casteras-Sourdis-2020}. In our setting,
we prove both directly from the equations for the segregated profile.}

		\subsection{The main results}
\revised{
A constant drift in the one-dimensional problem can arise from a change to
a moving frame. For related uses of transport terms in the study of rotating
waves, see \cite{Angenent-Fiedler-1988,Terracini-Verzini-Zilio-2019,Salort-Terracini-Verzini-Zilio-2025,Li-Terracini-2025,spadaro2026selection}.
We study the one-dimensional system}

	\begin{equation}\label{o-eq}
		\begin{cases}
			-u_i''-cu_i'=\mu u_i(1-u_i)-\beta u_i\displaystyle\sum_{j\ne i}\alpha_{ij}u_j,
			&\mbox{in } (-\uppi, \uppi),\\
			u_i(\pm\uppi)=0,
		\end{cases}
		\qquad i=1,2, \cdots,k,
	\end{equation}
	and the radial system
	\begin{equation}\label{o-system-radial}
		\begin{cases}
			- u_i''- \frac{n-1}{r} u_i' =\mu u_i(1-u_i)-\beta u_i\displaystyle\sum_{j\ne i}\alpha_{ij}u_j,
			&\mbox{in }  (0,R),\\
			u_i'(0)=0,\quad u_i(R)=0,
		\end{cases}
			\qquad i=1,2, \cdots,k,
	\end{equation}
	where $\mu, R,\beta>0$, $c\in\R$,  $\alpha_{ij}>0$ ($i\ne j$) are constants, $k\ge 2$ is the number of species, $n\ge 2$ is the space dimension. %We write $\bv u =(u_1,\dots,u_k)$. 
	
\revised{On each interval occupied by a single species in a segregated profile, the
density solves a logistic Dirichlet problem. For $a<b$ we write this problem as}
\begin{equation}\label{1d-1}
	-w''-cw'=\mu w(1-w)\quad\text{in }(a,b),\qquad w(a)=w(b)=0.
\end{equation}
\revised{We call the positive solution, when it exists, the logistic solution
on $(a,b)$. We begin with profiles in which each species occupies
one interval, in the order $1,\ldots,k$.} 
\begin{definition}\label{def-profile}
	A segregated profile is a vector $\bv U=(U_1, U_2, \cdots,U_k)$ of continuous
	functions on $[-\uppi,\uppi]$, together with points
	$-\uppi=s_0<s_1<\dots<s_k=\uppi$, such that for every $i$ the function $U_i$ is
	the logistic solution on $I_i=(s_{i-1},s_i)$, extended by zero outside $I_i$,
	and such that for $1\le i\le k-1$
	\begin{equation}\label{flux}
		\alpha_{i+1,i}\,p_i=\alpha_{i,i+1}\,q_i,\qquad
		p_i=-U_i'(s_i),\quad q_i=U_{i+1}'(s_i),
	\end{equation}
	the derivatives being taken inside the respective intervals.
	We call $s_1,s_2, \cdots,s_{k-1}$ the interface points, and we put $L_i=s_i-s_{i-1}$,
	$\vartheta_i=\alpha_{i+1,i}/\alpha_{i,i+1}$, $q_0=U_1'(-\uppi)>0$ and
	$p_k=-U_k'(\uppi)>0$.
\end{definition}

\revised{
We first construct a positive solution for a profile in which each species
occupies one interval. At an interface the segregated densities vanish
linearly. Balancing diffusion with competition therefore gives the scale
$\varepsilon=\beta^{-1/3}$ for both the distance from the interface and
the size of the two densities there. We construct a two-component
expansion at that scale and use the Taylor expansions of the segregated
densities to join it to the solutions away from the interface.

To estimate the inverse of the linearized operator, we take a weighted
difference of the two equations at each interface. The mutual competition
terms cancel, which permits passage to the limit across the interface.
Nondegeneracy of the segregated profile forces those limiting differences
to vanish. Scalar comparison estimates then give bounds on the individual
components. With these bounds, we use a contraction to correct the
approximate solution and a maximum principle to prove positivity.}

\begin{theorem}\label{main-theorem}
	Let $k\ge2$, $\mu>(c^2+k^2)/4$, $c\in\R$ and $\alpha_{ij}>0$ ($i\ne j$) be fixed,
	and let $\bv U$ be a segregated profile, as in Definition \ref{def-profile}.
	Then there are constants $\beta_0,C_1,C_2>0$ such that, for every $\beta\ge\beta_0$,  system \eqref{o-eq} has a solution
	$\bv u_\beta$ with all components strictly positive in $(-\uppi,\uppi)$, and
	\begin{equation}\label{global-convergence}
		C_1\beta^{-1/3}\le\|\bv u_\beta-\bv U\|_\infty\le C_2\beta^{-1/3}.
	\end{equation}
	Set $\varepsilon=\beta^{-1/3}$. For $1\le i\le k-1$ let $\bv a^{(i)}$ be the
	model profile of Proposition \ref{inner-profile} with
	$(\alpha_{12},\alpha_{21},p,q)=(\alpha_{i,i+1},\alpha_{i+1,i},p_i,q_i)$,
	and put $h=r-s_i$, $z=h/\varepsilon$. For every fixed $M'>0$, uniformly for
	$|h|\le M'\varepsilon\log\beta$,
	\begin{equation}\label{main-inner-expansion}
		u_{\beta,i}(r)=\varepsilon a_1^{(i)}(z)+O(\varepsilon^2+h^2),\qquad
		u_{\beta,i+1}(r)=\varepsilon a_2^{(i)}(z)+O(\varepsilon^2+h^2),
	\end{equation}
	while $|u_{\beta,j}(r)|\le Ce^{-1/(C\varepsilon)}$ for $j\notin\{i,i+1\}$, with $C>0$ depending on the fixed profile. 
	In particular,
	\begin{equation}\label{main-inner-limit}
		\varepsilon^{-1}(u_{\beta,i},u_{\beta,i+1})(s_i+\varepsilon z)\longrightarrow\bv a^{(i)}(z)
		\quad\text{in }C^2_{\mathrm{loc}}(\R).
	\end{equation}
	%
	%	The solution is unique in the neighborhood of the approximate solution
	%described in Proposition \ref{existence}; no global uniqueness for
	%	\eqref{o-eq} is asserted.
\end{theorem}
\begin{remark}
By Proposition \ref{outer-existence},	the assumption $\mu>(c^2+k^2)/4$ is optimal for the existence of a segregated profile defined in Definition \ref{def-profile}. 
\end{remark}

\revised{
The same construction allows a species to occupy several disjoint
intervals. We assign each interval a species label, require that all species
occur, and require different labels on consecutive intervals.}

\begin{definition}\label{def-labelled-profile}
	A  labelling is a surjective map $\sigma:\{1,2,\cdots,m\}\to\{1,2,\cdots,k\}$
	with $\sigma(l)\ne\sigma(l+1)$ for $1\le l\le m-1$. A  segregated profile for \eqref{o-eq} 
	with labelling $\sigma$ is given by points $-\uppi=s_0<s_1<\dots<s_m=\uppi$ and
	functions $V_l\in C^2(\bar I_l)$, $I_l=(s_{l-1},s_l)$, such that each $V_l$ is
	the logistic solution on $I_l$, and such that for $1\le l\le m-1$,
	\begin{equation}\label{flux-sigma}
		\alpha_{\sigma(l+1),\sigma(l)}\,p_l=\alpha_{\sigma(l),\sigma(l+1)}\,q_l,\qquad
		p_l=-V_l'(s_l)>0,\quad q_l=V_{l+1}'(s_l)>0.
	\end{equation}
	The profile of species $j$ is the function $U_j$ equal to $V_l$ on $I_l$
	whenever $\sigma(l)=j$, and to zero elsewhere. We write
	\begin{equation*}
			\Omega_j=\bigcup_{\sigma(l)=j}\bar I_l
	\end{equation*}
	for the set occupied by species $j$. Let $\Sigma_j$ be the set of interface points that
	are endpoints of intervals contained in $\Omega_j$. 
\end{definition}

\revised{Distinct labellings give distinct limits, so the solutions constructed from them are distinct for sufficiently large $\beta$.}  
\begin{theorem}\label{main-theorem-sigma}
	Let $m\ge k\ge2$, $\mu>(c^2+m^2)/4$, $c\in\R$ and $\alpha_{ij}>0$ ($i\ne j$) be fixed, let
	$\sigma$ be a labelling of $m$ intervals, and let
	$\bv U$ be the corresponding segregated profile, as stated in Definition \ref{def-labelled-profile}. For all sufficiently large
	$\beta$, system \eqref{o-eq} has a solution $\bv u_\beta$ with all components
	strictly positive in $(-\uppi,\uppi)$ 
\revised{satisfying \eqref{global-convergence}. At each $s_l$,
\eqref{main-inner-expansion} and \eqref{main-inner-limit} hold for
species $\sigma(l)$ and $\sigma(l+1)$ with the corresponding
coefficients and slopes. The remaining components are exponentially
small there, as in Theorem~\ref{main-theorem}.} 
	\revised{Let $m_0$ be the largest integer strictly less than $2\sqrt{\mu-c^2/4}$.
	For all sufficiently large $\beta$, system \eqref{o-eq} has at least} 
	\begin{equation*}
		\sum_{m=k}^{m_0} \sum_{i=0}^{k}(-1)^i \binom ki(k-i)(k-i-1)^{m-1}
	\end{equation*}
	distinct positive solutions.
\end{theorem}
\begin{remark}
\revised{
The count follows by inclusion--exclusion over the species labels omitted.
When a species occupies several intervals, its absence between two of
them requires a comparison estimate with data at both ends of the gap.
Lemma~\ref{two-sided-barrier} gives that estimate.}

\end{remark}

\revised{
For the radial problem, translation changes the coefficient $(n-1)/r$,
so the length of a shell alone does not determine its logistic solution.
We instead use the value at the center as the shooting parameter and
determine successive interfaces by matching the weighted slopes. We write
$j_{(n-2)/2,m}$ for the $m$th positive zero of the Bessel function
$J_{(n-2)/2}$.}

\begin{definition}\label{def-profile-radial}
	A labelling is a surjection $\sigma:\{1,2,\cdots,m\}\to\{1,2,\cdots,k\}$ with \revised{$\sigma(l)\ne\sigma(l+1)$ for $1\le l<m$. A segregated profile} for \eqref{o-system-radial} with this labeling consists of radii
	\begin{equation*}
			0=s_0<s_1<\cdots<s_m=R
	\end{equation*}
	and positive functions $V_l$ on their respective shells $I_l=(s_{l-1},s_l)$, satisfying
	\begin{equation*}
		-V_l''-\frac{n-1}{r} V_l' =\mu V_l(1-V_l),
	\end{equation*}
	with $V_1'(0)=0$, $V_1(s_1)=0$, and with zero values at both endpoints for $l\ge2$. \revised{At each interface $s_l$, require the weighted slopes to satisfy}
	\begin{equation}\label{flux-condi}
		\alpha_{\sigma(l+1),\sigma(l)}p_l
		=\alpha_{\sigma(l),\sigma(l+1)}q_\ell, \quad  	p_l=-V_l'(s_l)>0,\qquad q_l=V_{l+1}'(s_l)>0.
	\end{equation}
	The species profile $U_i$ equals $V_l$ on each shell labeled $i$, and is zero elsewhere. Write
	\begin{equation*}
			\Omega_i=\bigcup_{\sigma(l)=i}[s_{l-1},s_l],\qquad
		\Sigma_i=\partial\Omega_i\cap(0,R).
	\end{equation*}
	The first shell includes the center: $U_{\sigma(1)}(0)=V_1(0)>0$.
\end{definition}

\revised{For each radial profile, the next theorem gives a family of positive
solutions and describes its behavior near the interfaces.} 
	\begin{theorem}\label{thm-main-radial}
		Assume $	 \mu>j_{(n-2)/2,m}^2/R^2$ with $m\ge k\ge 2$, $\alpha_{ij}>0$ ($i\ne j$) fixed,   and fix a labelling $\sigma$ and its segregated profile $\bx U$, as defined in Definition  \ref{def-profile-radial}. Then there are constants $\beta_0,C_1,C_2>0$ such that, for every $\beta\ge\beta_0$, the radial system \eqref{o-system-radial} has a radial positive solution $\bx u_\beta$ satisfying
		\begin{equation*}
		%	0<u_{\beta,i}(r)<1\quad(0\le r<R),\qquad
			%
			C_1\beta^{-1/3} \le\|\bx u_\beta-\bx U\|_\infty\le C_2\beta^{-1/3}.
		\end{equation*}
		 Near each interface point $s=s_\ell$, let $\bv a^{(l)}$ be the
		 model profile of Proposition \ref{inner-profile} with
		 $(\alpha_{12},\alpha_{21},p,q)=(\alpha_{\sigma(l),\sigma(l+1)},\alpha_{\sigma(l+1),\sigma(l)},p_l,q_l)$,
		 and put $h=r-s_\ell$, $z=h/\varepsilon$ with $\eps=\beta^{-\frac13}$. Then
		\begin{equation*}
			\eps^{-1}(u_{\beta,\sigma(l)},u_{\beta,\sigma(l+1)})(s+\eps z)\longrightarrow\bx a^{(l)}(z)
			\quad\text{in }C^2_{\mathrm{loc}}(\R).
		\end{equation*}
		For any fixed $M'>0$, uniformly for $|h|\le M'\eps\log\beta$,
		\begin{equation*}
			(u_{\beta,\sigma(l)},u_{\beta,\sigma(l+1)})(s+h)
			=\eps\bx a^{(l)}(h/\eps)+O(\eps^2+h^2).
		\end{equation*}
		Every other component is bounded there by $Ce^{-c_*/\eps}$, with $c_*>0$ depending on the fixed profile. 
		Set $	M_\mu=\#\{m\ge1:j_{(n-2)/2,m}<R\sqrt\mu\}$.
		For all sufficiently large $\beta$, the number of distinct positive radial solutions is at least
		\begin{equation*}
			\sum_{m=k}^{M_\mu} \sum_{i=0}^k(-1)^i\binom ki \revised{(k-i)(k-i-1)^{m-1}}.
		\end{equation*}
	\end{theorem}
	\begin{remark}
		The assumption $\mu>j_{(n-2)/2,m}^2/R^2$ is optimal by Proposition \ref{profile-radial}. 
	\end{remark}

\revised{We finally consider a species $j$ which experiences no competition from
species $l$. On an interval or shell occupied by $l$, the other species
can still be suppressed by competition with $l$. This leads to the
following obstruction.}

\begin{proposition}\label{prop-zero}
\revised{Let $j\ne l$, let $\alpha_{jl}=0$, and assume
$\alpha_{hl}>0$ for every $h\notin\{j,l\}$. Assume all other
competition coefficients are nonnegative. Let $\bv U$ be a segregated
profile for the interval or the ball, with the endpoint conditions and
flux balances of the corresponding definition for these coefficients.
Then no sequence of positive solutions with $\beta\to\infty$ can
converge uniformly to $\bv U$.}
\end{proposition}

	\subsection{Notation.}
	Throughout this paper, we adopt the following notation.
	\begin{enumerate}
		\item  $n\ge2$ is the space dimension, $k\ge2$ is the number of species, and $m\ge k$ is the number of occupied intervals or radial shells. Let $\upnu=(n-2)/2$ and let $j_{\upnu,m}$ be the $m$th positive zero of the Bessel function $J_\upnu$. 
		\item  \revised{$\varepsilon=\beta^{-1/3}$, $\rho=\varepsilon\log\beta$ and
$\varrho(r)=r^{n-1}$. For vector-valued functions,
$\|\bv v\|_\infty=\max_j\|v_j\|_\infty$.}
		\item The small positive constant $\zeta$ is fixed, satisfying $0<\zeta\le \min\{\zeta_0, \zeta_1,\zeta_2\}$, where $\zeta_0,
		\zeta_1$, $\zeta_2$  are constant given by comparison theorem in Lemmas \ref{scalar-barriers-1},~\ref{two-sided-barrier},~\ref{gap-left-out-barrier}, respectively. 
		\item  Let $\chi\in C_c^\infty(\R)$ satisfy $0\le\chi\le1$, $\chi=1$ on $[-M_0,M_0]$, 
		$\chi=0$ outside $[-M,M]$, and put $\theta_i(r)=\chi((r-s_i)/\rho)$ for
		$1\le i\le k-1$. 
		\item We say that a function is
		{\it rapidly decaying} at an end if, for every $\tau>0$ and every fixed
		derivative order $j$, its $j$th derivative is $O(e^{-\tau|z|})$ at that end.  Each occurrence of $\mathcal E_\pm$ denotes a possibly different
		function that is  rapidly decaying at the corresponding end $\pm \infty$, together with
		its derivatives. 
		\item In our construction process, we
		choose $\tau>\zeta$ and $0<M_0<M$ such that $M_0(\tau-\zeta)>1$.
		%	\item $\sigma_\beta=\sum_lu_l^{\ap}$.
		%	\item $\varsigma\in (0,1)$ is a fixed constant, and $\varsigma_*=3+\varsigma$.
	\end{enumerate}	
	
\revised{
Section~\ref{1-d-m-eq-k} proves the construction for one interval per
species. Section~\ref{1-d-multiplicity} treats repeated labels, and
Section~\ref{radial-sys} proves the radial results.
Section~\ref{vanish-rate} proves Proposition~\ref{prop-zero}.
The appendices give the scalar logistic and comparison results used in
these proofs, including the shooting argument for the radial profiles.}

\section{\revised{The interval problem with one interval per species}} \label{1-d-m-eq-k}

\subsection{Existence and nondegeneracy of segregated profiles}
By Lemma \ref{eigenvalue-n}, the Dirichlet eigenvalues of $-\frac{\mathrm{d}^2}{\mathrm{d}x^2} - c \frac{\mathrm{d}}{\mathrm{d}x} $ on $(a,b)$ are given by 
\begin{equation*}
\lambda_j(a,b)=\big(\frac{j \uppi}{b-a}\big)^2+\frac{c^2}{4}, \qquad j=1,2, \cdots, 
\end{equation*}
with eigenfunctions $\mathrm{e}^{-\frac{c}{2}x}\sin \frac{j\uppi(x-a)}{b-a}$.  Furthermore, in weighted space $L^2\big((a,b), \mathrm{e}^{cx} \mathrm{d}x \big)$, the operator is self-adjoint. All integrations by parts are carried out in this space.

\begin{proposition}\label{outer-existence}
	A segregated profile defined in Definition \ref{def-profile}  exists if and only if $\mu-c^2/4>k^2/4$. \revised{For the fixed ordering in Definition~\ref{def-profile}, it is unique.}
\end{proposition}
\begin{proof}
\revised{
Existence of a positive logistic solution on each interval requires
$\mu>c^2/4$ by Proposition~\ref{exist-unique}. The same proposition
then gives $L_i>L_*=\uppi/\sqrt{\mu-c^2/4}$. Summing over the
intervals yields $2\uppi>kL_*$, which is the claimed condition.
For the converse, assume $\mu>(c^2+k^2)/4$.}

By uniqueness of the logistic solution and
	translation invariance, $U_i$ is the translate of the solution $w_{L_i}$
	of Lemma \ref{slope-monotonicity}. Hence $p_i=A(L_i)$ and $q_i=B(L_{i+1})$. \revised{Thus the lengths determine the interval endpoints and the logistic
solutions. We seek lengths $L_i>L_*$ that sum to $2\uppi$ and satisfy
the following form of the flux balance.}
	\begin{equation}\label{length-recursion}
		B(L_{i+1})=\vartheta_iA(L_i),\qquad 1\le i\le k-1.
	\end{equation}
	By Lemma \ref{slope-monotonicity}, the mapping $B:(L_*,\infty)\to (0,B_\infty)$ is
	a continuous increasing bijection, with
	$B_\infty=\sup B\in(0,\infty]$, and $A$ is continuous and increasing with
	$A(L)\to0$ as $L\downarrow L_*$. Now we choose the first length $L_1$
	freely and compute the others from \eqref{length-recursion}. Given
	$L_1>L_*$, \eqref{length-recursion} determines $L_2,L_3,\dots$ one after the
	other, as long as $\vartheta_iA(L_i)<B_\infty$. Let $D_i$ be the set of
	$L_1$ for which $L_1,\cdots,L_i$ are defined. Since $A$ and $B$ are increasing,
\revised{each $L_i$ is a continuous, strictly increasing function of $L_1$
where it is defined. The condition $\vartheta_iA(L_i)<B_\infty$
restricts $D_i$ to an interval with left endpoint $L_*$. Inductively,
as $L_1\downarrow L_*$, each $L_i\downarrow L_*$ because $A(L_i)$
and $B(L_i)$ tend to zero there.}
	
 It follows that $\Upsilon(L_1):=\sum_iL_i$ is continuous and
	strictly increasing on $D_k=(L_*,L^\sharp)$, and that $\Upsilon\to kL_*$ as
	$L_1\downarrow L_*$. We claim that $\Upsilon\to\infty$ as $L_1\uparrow
	L^\sharp$. \revised{If $L^\sharp=\infty$, this follows from
$\Upsilon\ge L_1$. Otherwise, let $j+1$ be the first index with
$\sup D_{j+1}=L^\sharp$. The length $L_j$ is defined and continuous
at $L^\sharp$. At this endpoint the condition for continuing the
recursion reaches equality, so}
	$\revised{\vartheta_jA(L_j(L_1))\uparrow B_\infty}$ as $L_1\uparrow L^\sharp$, so that
	\revised{\revised{$L_{j+1}\to\infty$. Hence} $\Upsilon$ tends to infinity as $L_1\uparrow L^\sharp$, and $\Upsilon(L_1)=2\uppi$ has exactly one}
	solution when $kL_*<2\uppi$. Since the lengths determine the profile, the
	profile exists and is unique.
\end{proof}

\begin{lemma}[Nondegeneracy]\label{nondegeneracy}
	Let $\Phi_i\in C^2(\bar I_i)$, $1\le i\le k$, satisfy $L_{U_i}\Phi_i=0$ in
	$I_i$, $\Phi_1(-\uppi)=\Phi_k(\uppi)=0$, and, for $1\le i\le k-1$,
	\begin{equation}\label{kernel-conditions}
		\alpha_{i,i+1}\Phi_{i+1}(s_i)=-\alpha_{i+1,i}\Phi_i(s_i),\qquad
		\alpha_{i,i+1}\Phi_{i+1}'(s_i)=-\alpha_{i+1,i}\Phi_i'(s_i).
	\end{equation}
	Then $\Phi_i=0$ for every $i$.
\end{lemma}
\begin{proof}
	Put
	\begin{equation*}
			w_i=\alpha_{i,i+1}\Phi_{i+1}(s_i)
		=-\alpha_{i+1,i}\Phi_i(s_i),\qquad 1\le i\le k-1.
	\end{equation*}
	By uniqueness for the Dirichlet problems and \eqref{Epm}, one sees that
	\begin{equation}\label{Phi-i-1}
			\Phi_i=\frac{w_{i-1}}{\alpha_{i-1,i}}E_i^--\frac{w_i}{\alpha_{i+1,i}}E_i^+, 
	\end{equation}
	where the first term is absent for $i=1$ and the second for $i=k$. Inserting
	this into the second condition in \eqref{kernel-conditions} gives
	$\mathsf Mw=0$, where $\mathsf M$ is the tridiagonal $(k-1)\times(k-1)$ matrix
	with entries
	\begin{equation*}
			\mathsf M_{ii}=(E_{i+1}^-)'(s_i)-(E_i^+)'(s_i),\quad
		\mathsf M_{i,i+1}=-\frac{\alpha_{i,i+1}}{\alpha_{i+2,i+1}}(E_{i+1}^+)'(s_i),\quad
		\mathsf M_{i,i-1}=\frac{\alpha_{i+1,i}}{\alpha_{i-1,i}}(E_i^-)'(s_i).
	\end{equation*}
	By Lemma \ref{transversality} the diagonal entries are positive, and by
	\eqref{Epm-hopf} the off-diagonal entries are negative.
	
	Let
	$\Theta_i=-U_i'$. It satisfies
	$L_{U_i}\Theta_i=0$, $\Theta_i(s_{i-1})=-q_{i-1}$ and $\Theta_i(s_i)=p_i$. By
	\eqref{flux} and the relation $U_i''=-cU_i'$ at the zeros, it satisfies both
	conditions in \eqref{kernel-conditions}, with $w_i$ replaced by $-x_i$, where
	\begin{equation*}
		x_i=\alpha_{i,i+1}q_i=\alpha_{i+1,i}p_i>0.
	\end{equation*}
	While we notice that $\Theta_1(-\uppi)=-q_0$ and
	$\Theta_k(\uppi)=p_k$. Let $\Phi[-x]$ denote the functions defined by \eqref{Phi-i-1} with $w=-x$. Then $\Phi_i[-x]=\Theta_i$ for $2\le
	i\le k-1$, while $\Phi_1[-x]=\Theta_1+q_0E_1^-$ and
	$\Phi_k[-x]=\Theta_k-p_kE_k^+$. Since $\Theta$ satisfies the derivative
	conditions exactly, only these two corrections contribute to $\mathsf M(-x)$.
	Therefore $(\mathsf Mx)_i=0$ for $2\le i\le k-2$, and by \eqref{Epm-hopf}
	\begin{equation*}
			(\mathsf Mx)_1=-\alpha_{21}q_0(E_1^-)'(s_1)>0,\qquad
		(\mathsf Mx)_{k-1}= \alpha_{k-1,k}p_k(E_k^+)'(s_{k-1})>0,
	\end{equation*}
	and the two terms add up in the single row when $k=2$.
	\revised{To see why this implies invertibility, let
$\mathsf D=\operatorname{diag}(x_1,\ldots,x_{k-1})$ and
$\mathsf B=\mathsf D^{-1}\mathsf M\mathsf D$. Its row sums are
nonnegative and the first and last are strictly positive. If
$\mathsf Bv=0$ and $|v_i|=\|v\|_\infty>0$, the equation in row $i$
and the signs of the off-diagonal entries force equality in the
triangle inequality. That row must have zero sum and every neighboring
component must have the same absolute value. Repeating this argument
along the consecutive rows reaches a row with positive sum, a
contradiction. Thus $\mathsf M$ has trivial kernel. In particular
$w=0$, and \eqref{Phi-i-1} gives $\Phi_i=0$ for every $i$.}
\end{proof}

\subsection{Model profile near an interface point}\label{sec-inner}
\revised{
We first work near $s_1$. At this interface the limiting profiles of
species $1$ and $2$ have nonzero one-sided slopes. The remaining species
are absent on a fixed neighborhood. We use the variable
$z=(r-s_1)/\varepsilon$ and seek an approximation to the first two
densities of the form}

\begin{equation}\label{model-ansatz}
	\eps\,\bv a(z)+\eps^2\bv P(z)+\eps^3\bv Q(z), \qquad \eps =\beta^{-\frac13}. 
\end{equation}
\revised{The leading term solves the two-species equations obtained by keeping
diffusion and competition at order $\varepsilon^{-1}$. Its slopes at
infinity must agree with those of the segregated profile.} 
\begin{proposition}\label{inner-profile}
	Put $K=\alpha_{21} p=\alpha_{12} q>0$. There is a unique positive profile
	$\bv a=(a_1,a_2)$ satisfying
	\begin{equation}\label{a-eq}
		\begin{aligned}
			\begin{cases}
									a_1''=\alpha_{12} a_1a_2,\quad & a_2''=\alpha_{21} a_1a_2,
				\\
					a_1(z)\sim-pz  \ \mbox{at } -\infty, \quad& \revised{a_2(z)\to 0} \ \mbox{at } -\infty,
				\\
				\revised{a_1(z)\to 0} \  \mbox{at } +\infty, 	\quad  & a_2(z)\sim qz \ \mbox{at } +\infty,
			\end{cases}
		\end{aligned}
	\end{equation}
and
	\begin{equation}\label{center}
		\alpha_{12} a_2(z)-\alpha_{21} a_1(z)=Kz.
	\end{equation}
	It satisfies $a_1'<0<a_2'$.
\end{proposition}
\begin{proof}
	By \cite[Lemma 4.2]{2005CTV}, there is a unique solution of
	$\eta''=\eta(\eta-y)$ on $\R$ with $\eta(y)>y^+$.
	\revised{The identity $\eta(y)=y+\eta(-y)$ relates the two ends:
$\eta(y)\to0$ as $y\to-\infty$ and $\eta(y)-y\to0$ as
$y\to+\infty$. Moreover, $\eta''>0$ and $\eta'\ge0$ by convexity
and monotonicity, so $\eta'>0$ at every finite point.}
	For $\kappa=K^{1/3}$, set
	\begin{equation*}
			a_1(z)=\frac{\kappa^2}{\alpha_{21}}\eta(-\kappa z),\qquad
		a_2(z)=\frac{\kappa^2}{\alpha_{12}}\eta(\kappa z).
	\end{equation*}
	A direct computation, using $\eta(-y)=\eta(y)-y$ and $\kappa^3=K$, shows that
	these functions satisfy \eqref{a-eq}, have the required growth, and satisfy
	\eqref{center}.
	
	For uniqueness, let $\bv a$ be any positive solution with the stated
	asymptotics. By \eqref{a-eq}, $(\alpha_{12} a_2-\alpha_{21} a_1)''=0$, and the
	slopes at $\pm\infty$, together with $\alpha_{21} p=\alpha_{12} q=K$, force
	$\alpha_{12} a_2-\alpha_{21} a_1=Kz+d$ for some $d\in\R$. Replacing $\bv a(z)$
	by $\bv a(z-d/K)$ removes $d$ without changing the slopes, so \eqref{center}
	fixes $d=0$. Then $\eta(y):=(\alpha_{12}/\kappa^2)a_2(y/\kappa)$ satisfies
	$\eta-y=(\alpha_{21}/\kappa^2)a_1(y/\kappa)>0$, hence $\eta>y^+$ and
	$\eta''=\eta(\eta-y)$, and the uniqueness part of the cited lemma gives the
	claim.
\end{proof}
\begin{remark}
	\revised{At $s_i$ the same profile is obtained by using the coefficients
$\alpha_{i,i+1}$ and $\alpha_{i+1,i}$ and the slopes $p_i$ and $q_i$.}
\end{remark}

 The lemma below compares the local expansion with the segregated profile.
\begin{lemma}\label{decay}
	For every $\tau>0$ and every fixed integer $j\ge0$, 
	\begin{equation*}
		\begin{aligned}
			\big|\partial_z^j(a_1+pz)\big|+\big|\partial_z^ja_2\big|
		&\le C_{\tau,j}e^{-\tau|z|}, \quad z\le-1,\\
		\big|\partial_z^ja_1\big|+\big|\partial_z^j(a_2-qz)\big|
		&\le C_{\tau,j}e^{-\tau z}, \quad \quad z\ge1.
		\end{aligned}
	\end{equation*}
\end{lemma}
\begin{proof}
	At $+\infty$, $a_1$ tends to zero and solves the linear equation
	$-a_1''+\alpha_{12} a_2a_1=0$, whose coefficient $\alpha_{12} a_2(z)\sim Kz$
	grows without bound. Given $\tau>0$, choose $R$ so large that
	$\alpha_{12}a_2\ge\tau^2$ on $[R,\infty)$. Then the function $Ce^{-\tau(z-R)}$
	is a supersolution there, and with $C\ge a_1(R)$ the comparison principle
	gives $a_1\le Ce^{-\tau(z-R)}$. Since $\tau$ is arbitrary, $a_1$ decays faster
	than any exponential. The same bound for the derivatives follows from the
	equation and from interpolation between the function and its second
	derivative. Differentiating may produce polynomial factors in $z$, but these
	do not affect the conclusion: it suffices to apply the estimate with a
	slightly larger rate. Identity \eqref{center} gives the estimate for $a_2-qz$. The
	proof at $-\infty$ is the same, with the components interchanged.
\end{proof}

\revised{
To match the Taylor terms through degree three, we solve for two
corrections to the leading profile. Their equations use the linearization
of \eqref{a-eq} at $\bv a$. With
$\bv b=(\alpha_{12},\alpha_{21})^T$, this operator is}

\begin{equation}\label{inner-operator}
	\M\bv X=-\bv X''+\bv b\,(a_2X_1+a_1X_2), \quad \bv X=(X_1,X_2).
\end{equation}
That is,
\begin{equation*}
	(\M\bv X)_1=-X_1''+\alpha_{12}(a_2X_1+a_1X_2),\qquad
	(\M\bv X)_2=-X_2''+\alpha_{21}(a_2X_1+a_1X_2).
\end{equation*}
We further define 
\begin{equation*}%\label{difference-def}
	\Dif{\bv X}=X_1/\alpha_{12}-X_2/\alpha_{21}.
\end{equation*}
Then 
\begin{equation*}
	\Dif{\M\bv X}= 	\frac{(\M\bv X)_1}{\alpha_{12}}-\frac{(\M\bv X)_2}{\alpha_{21}}
	=-\frac{X_1''}{\alpha_{12}}+\frac{X_2''}{\alpha_{21}}=-\Dif{\bv X}''.
\end{equation*}
Thus, if $\M\bv X$ is known, then $\Dif{\bv X}$ is determined up to an affine
function. For the profile $\bv a$,  identity \eqref{center} and
$K=\alpha_{21}p$ give
\begin{equation}
	\label{Da}
	\M\bv
	a=\bv a'', \quad  \M\bv a'=0, \quad 
	\Dif{\bv a}=\frac{\alpha_{21} a_1-\alpha_{12} a_2}{\alpha_{12}\alpha_{21}}=-\frac{pz}{\alpha_{12}}.
\end{equation}

The next lemma is the basic solvability result for $\M$, which is used to construct the corrections.

\begin{lemma}\label{inverse}
	Let $\bv H=(H_1,H_2)$ be smooth, with $H_2$ rapidly decaying at $-\infty$ and
	$H_1$ rapidly decaying at $+\infty$. Suppose a polynomial $D$ satisfies
	\begin{equation}\label{D-compatible}
		D''=-H_1/\alpha_{12}+H_2/\alpha_{21}.
	\end{equation}
	There is a unique smooth solution $\bv X$ of $\M\bv X=\bv H$ with
	$\Dif{\bv X}=D$ and
	\begin{equation}\label{Xleft-right}
		\begin{aligned}
			&	X_1=\alpha_{12} D+\mathcal E_-,\quad X_2=\mathcal E_-   \qquad \quad (z\to-\infty),
			\\
		&	X_1=\mathcal E_+,\qquad X_2=-\alpha_{21} D+\mathcal E_+,  \quad(z\to+\infty),
		\end{aligned}
	\end{equation}
	where each occurrence of $\mathcal E_\pm$ denotes a possibly different
	function that is rapidly decaying at the corresponding end, together with
	its derivatives.
\end{lemma}
\begin{proof}
	Set $S=\alpha_{21} a_1+\alpha_{12} a_2$ and
	$\Kop=-\partial_z^2+S$. Then $S>0$, $S\sim K|z|$ at infinity, and
	$\inf_\R S>0$. Choose a smooth function $\chi_+$ equal to zero on
	$(-\infty,-1]$ and to one on $[1,\infty)$, and define
	\begin{equation}\label{g-equation}
		g=H_2/\alpha_{21}-\alpha_{12} a_2D+\Kop(\chi_+D).
	\end{equation}
	At the left end, we have $g=H_2/\alpha_{21}-\alpha_{12}a_2D$;
	at the right end, \eqref{D-compatible} gives $g=H_1/\alpha_{12}+\alpha_{21}a_1D$.
	In both cases, $g$ and its derivatives
	are rapidly decaying.
\revised{We use the bilinear form $\int(\psi'\eta'+S\psi\eta)$ on the Hilbert space}
	\begin{equation*}
			\mathcal H=\{\psi\in H^1(\R):\int_\R S\psi^2<\infty\},
	\end{equation*}
	\revised{with squared norm $\int (|\psi'|^2+S|\psi|^2)$. Since $\inf S>0$, this
norm controls the $L^2$ norm and $g\in L^2(\R)$ defines a bounded linear
functional. Lax--Milgram gives a unique $\psi\in\mathcal H$ with
$\Kop\psi=g$. The equation gives smoothness, and $\psi\in H^1(\R)$
implies that $\psi$ tends to zero at both ends. Fix $\tau>0$ and choose
$R$ so that $S\ge\tau^2+1$ for $|z|\ge R$. There
$\Kop(e^{-\tau|z|})\ge e^{-\tau|z|}$. A large multiple dominates
$|g|$ under $\Kop$ and dominates $|\psi|$ at $z=\pm R$. Comparison on
the two half-lines then bounds $|\psi|$ by that multiple. The equation
and its derivatives give the same conclusion for each derivative, using
a larger rate for $\psi$ to absorb polynomial factors from $S$.}
	
Now set
	\begin{equation*}%\label{Xconstruction}
		X_1=\alpha_{12}\big((1-\chi_+)D+\psi\big),\qquad
		X_2=\alpha_{21}\big(\psi-\chi_+D\big).
	\end{equation*}
	\revised{The definition gives $\Dif{\bv X}=D$. Equation \eqref{g-equation}
gives $(\M\bv X)_2=H_2$, and \eqref{D-compatible} then gives
$(\M\bv X)_1=H_1$.} The asymptotic behavior \eqref{Xleft-right}
	follows from the definition, since $\psi$ is rapidly decaying.
	
If two solutions existed, their difference would
	satisfy $\Dif{\bv X}=0$, so it is of the form $(\alpha_{12} w,\alpha_{21}
	w)$, with $\Kop w=0$ and $w$ decaying. Testing by $w$ gives
	$\int(|w'|^2+Sw^2)=0$, so $w=0$.  Thus the solution is unique.  %In particular the solution does not depend on
%	the auxiliary choice of $\chi_+$.
\end{proof}
We can now construct the corrections in \eqref{model-ansatz}.  Recall that $\bv b=(\alpha_{12},\alpha_{21})^T$.

\begin{proposition}\label{correctors}
	Define
	\begin{equation*}
		\bv P(z)=-\frac{c}{10}\bigl(z^2\bv a'(z)+4z\bv a(z)\bigr).
	\end{equation*}
Then
	\begin{equation}\label{P-equation}
		\begin{aligned}
			\begin{cases}
					\M\bv P =c\bv a', \quad  
					\Dif{\bv P}=\frac{cp}{2b_1 }z^2,
					\\
					(P_1, P_2)= (cpz^2/2 +\mathcal E_-, \mathcal E_-) \  \mbox{at } -\infty,
			\\
				(P_1, P_2)=(	\mathcal E_+, -cqz^2/2+\mathcal E_+)
			\ \mbox{at } +\infty.
			\end{cases}
		\end{aligned}
	\end{equation}
		Moreover, there is a unique smooth $\bv Q$  such that
		\begin{equation}\label{Q-equation}
		\begin{aligned}
			\begin{cases}
					\M\bv Q=c\bv P'+\mu\bv a-\bv b\,P_1P_2, \quad
					 \Dif{\bv Q}=\frac{(\mu-c^2)p}{6 b_1 }z^3,
					\\
					(Q_1, Q_2)= \big( 	(\mu -c^2)pz^3/6 +\mathcal E_- , \mathcal E_-\big)
					\  \mbox{at } -\infty, \\
						(Q_1, Q_2)= \big(
						\mathcal E_+ ,  -	(\mu -c^2)qz^3/6+\mathcal E_+
					\big) \ \mbox{at } +\infty.
			\end{cases}
		\end{aligned}
		\end{equation}
 Here each occurrence of $\mathcal E_\pm$ denotes a possibly different
	function that is rapidly decaying at the corresponding end, together with
	its derivatives.
\end{proposition}
\begin{proof}
First, for a smooth function $f(z)$, straightforward calculation yields
\begin{equation*}
		\M(f\bv X)=f\M\bv X-2f'\bv X'-f''\bv X.
\end{equation*}
	Consequently, we get from \eqref{Da} that
	\begin{equation*}
			\M(z^2\bv a')=-4z\bv a''-2\bv a',\qquad \M(z\bv a)=z\bv a''-2\bv a'.
	\end{equation*}
	Hence $\M(z^2\bv a'+4z\bv a)=-10\bv a'$, which proves $\M\bv P=c\bv a'$.
	Identity \eqref{Da} gives $\Dif{\bv P}$, and Lemma \ref{decay} gives the
	behavior of $\bv P$ at infinity.
	
Put $\bv H=c\bv P'+\mu\bv a-\bv b\,P_1P_2$. Its
	second component decays rapidly at the left end, and its first component at
	the right end, because there they only involve the decaying components of
	$\bv a$ and $\bv P$. Moreover,
	\begin{equation*}
			\frac{H_1}{\alpha_{12}}-\frac{H_2}{\alpha_{21}}
		=c\Dif{\bv P}'+\mu\Dif{\bv a}=\frac{(c^2-\mu)p}{b_1}z.
	\end{equation*}
	Thus $D=(\mu-c^2)pz^3/(6 b_1)$ satisfies \eqref{D-compatible}, and Lemma
	\ref{inverse} gives $\bv Q$, uniquely, with the behavior \eqref{Q-equation}.
%	Since $D$ has no constant or linear part, neither does the polynomial
	%behavior of $\bv Q$.
\end{proof}

\subsection{The approximate solution}
For $1\le i\le k-1$ let $\bv a^{(i)}$, $\bv P^{(i)}$, $\bv Q^{(i)}$ be the
profile and the corrections with
$(\alpha_{12},\alpha_{21},p,q)=(\alpha_{i,i+1},\alpha_{i+1,i},p_i,q_i)$ in Subsection \ref{sec-inner}.
%All constants below are independent of $\beta$. 
%
Fix $1\le i\le k-1$, set $h=r-s_i$, $z=h/\eps$, and define the local expansion
\begin{equation}\label{inner}
	\bv u^{(i)}=\eps\,\bv a^{(i)}(z)+\eps^2\bv P^{(i)}(z)+\eps^3\bv Q^{(i)}(z), \quad \eps =\beta^{-\frac13}. 
\end{equation}
Its first component approximates $u_i$ and its second $u_{i+1}$. On the other
hand, the segregated profiles have the Taylor expansions
\begin{equation}\label{Vtaylor}
	\begin{aligned}
		U_i(s_i+h)&=-p_ih+\frac{cp_i}{2}h^2+\frac{(\mu-c^2)p_i}{6}h^3+O(h^4),
		\ 
		U_{i+1}(s_i+h)&=q_ih-\frac{cq_i}{2}h^2+\frac{(c^2-\mu)q_i}{6}h^3+O(h^4).
	\end{aligned}
\end{equation} 
The remainders, differentiated $j$ times, are
$O(|h|^{4-j})$ for $0\le j\le3$, on the side of $s_i$ where the respective
profile is positive.

\begin{proposition}\label{matching}
	For $0\le j\le3$, every fixed $\tau>0$, and $h<0$ sufficiently small,
	\begin{align}
		\big|\partial_r^j(u^{(i)}_1-U_i)\big|
		&\le C\big(|h|^{4-j}+\eps^{1-j}e^{-\tau|h|/\eps}\big),\label{match-left}\\
		\big|\partial_r^j u^{(i)}_2\big|
		&\le C\eps^{1-j}e^{-\tau|h|/\eps}.\label{weak-left}
	\end{align}
	For $h>0$, the corresponding estimates hold for $u^{(i)}_2-U_{i+1}$ and
	$u^{(i)}_1$. In particular, let $\ell=\log\beta$, $\rho=\eps\ell$, and fix
	$0<M_0<M$. If $\tau M_0\ge1$, then on the two intervals where $M_0 \rho\le|h|\le
	M\rho$, which we call the interpolation regions, the component that is
	present satisfies
	\begin{equation}\label{matching-final}
		\big|\partial_r^j(\text{local expansion}-\text{segregated profile})\big|
		\le C\eps^{4-j}\ell^{4-j},\qquad 0\le j\le3.
	\end{equation}
\end{proposition}
\begin{proof}
	On the left, by \eqref{a-eq}, \eqref{P-equation}, and \eqref{Q-equation},  the polynomial parts of the three terms in
	\eqref{inner} are exactly $-p_ih+\frac{cp_i}{2}h^2+\frac{(\mu-c^2)p_i}{6}h^3$, and
	by Lemma \ref{decay} the remaining parts are rapidly decaying in $z$. Hence
	\begin{equation*}
		\partial_r^j \Big[ u^{(i)}_1+p_ih-\frac{cp_i}{2}h^2
	-\frac{(\mu-c^2)p_i}{6}h^3 \Big]
	=O(\eps^{1-j}e^{-\tau|h|/\eps}),
	\end{equation*}
	where the factor $\eps^{1-j}$ accounts for the overall factor $\eps$ in
	\eqref{inner} and for $\partial_r=\eps^{-1}\partial_z$. Combining this with
	\eqref{Vtaylor} proves \eqref{match-left}, and \eqref{weak-left} follows in
	the same way, since the second component has no polynomial part on the left.
	\revised{On the right, the second component has the polynomial part and the
first is rapidly decaying, which gives the corresponding estimates.
In the interpolation regions,}
	$e^{-\tau|h|/\eps}\le\revised{\beta^{-\tau M_0}}\le\beta^{-1}=\eps^3$, while
	$|h|^{4-j}\le(M\rho)^{4-j}$, which gives \eqref{matching-final}.
\end{proof}

Let $\chi\in C_c^\infty(\R)$ satisfy $0\le\chi\le1$, $\chi=1$ on $[-M_0,M_0]$ and
$\chi=0$ outside $[-M,M]$, and put $\theta_i(r)=\chi((r-s_i)/\rho)$ for
$1\le i\le k-1$. Thus $\theta_i$ equals one on the region $|r-s_i|\le M_0\rho$,
where we use the local expansion, and zero outside $|r-s_i|\le M\rho$. For large $\beta$,  the supports of the $\theta_i$ are
disjoint and contained in $(-\uppi,\uppi)$. Extend each $U_l$ by zero outside
$I_l$, and define, for $1\le l\le k$, 
\begin{equation}\label{gluing}
	u_l^{\ap}=(1-\theta_{l-1}-\theta_l)U_l+\theta_{l-1}u^{(l-1)}_2+\theta_lu^{(l)}_1,
\end{equation}
with the convention that the terms involving $\theta_0$ or $\theta_k$ are
absent. In particular, $u_l^{\ap}$ vanishes
at distance $M\rho$ from $\bar I_l$, so that  $u_l^{\ap}u_j^{\ap}\equiv0$
whenever $|l-j|\ge2$.

\revised{To apply scalar comparison later, we need lower bounds both for
the component occupying a given interval and for the total density.}
\begin{lemma}\label{profile-geometry}
	For $1\le j,l \le k$ let
	$\Sigma_j=\{s_{j-1},s_j\}\cap(-\uppi,\uppi)$ be the set of interface points of $I_j$,
	and put $\sigma_\beta=\sum_l u_l^{\ap}$, where $u_l^{\ap}$ is defined by \eqref{gluing}. Then
	there are constants $\xi,C>0$, %determined by the fixed profiles, 
	such that,
	for all sufficiently large $\beta$, the following holds.
	\begin{enumerate}
		\item[(a)] $\|\bv u^{\ap}-\bv U\|_\infty\le C\eps$ and
		$\|\bv u^{\ap}\|_\infty\le C$.
		\item[(b)] For each $l$ and $r\in\bar I_l$,
		$u_l^{\ap}(r)\ge \xi \min\{\eps+\dist(r,\Sigma_l),\,r+\uppi,\,\uppi-r\}$.
		\item[(c)] For each $l$, each fixed $\tau>0$ and $r\notin I_l$,
		$|u_l^{\ap}(r)|\le C_\tau\eps e^{-\tau\dist(r,I_l)/\eps}$; moreover
		$u_l^{\ap}(r)=0$ if $\dist(r,I_l)\ge M\rho$, and $u_l^{\ap}\ge-C\eps^2$
		everywhere.
		\item[(d)] $\sigma_\beta\ge\xi\min\{r+\uppi,\uppi-r,\eps\}$, and
		$|u_l^{\ap}|\le C\sigma_\beta$ for each $l$.
		\item[(e)] For each $l$ and $r\notin I_l$,
		$\sum_{j\ne l}\alpha_{lj}u_j^{\ap}(r)\ge\xi\min\{\eps,r+\uppi,\uppi-r\}$ and
		$|u_l^{\ap}(r)|\le C\sum_{j\ne l}u_j^{\ap}(r)$.
	\end{enumerate}
\end{lemma}
\begin{proof}
	Away from the sets $\{|r-s_i|\le M\rho\}$ one has $\bv u^{\ap}=\bv U$, and all
	claims follow from the positivity of the $U_l$ on their intervals and from
	Hopf's lemma at $\pm\uppi$. Near $s_i$ only the two model components are
	nonzero; we write $h=r-s_i$ and drop the superscript $(i)$.
For $z\ge0$, positivity and the linear
	growth of $a_2$ give $a_2(z)\ge\xi_0(1+z)$, whereas $|P_2(z)|\le C(1+z^2)$ and
	$|Q_2(z)|\le C(1+z^3)$. On $0\le z\le M\log\beta$ the corrections are
	relatively small after multiplication by $\eps^2$ and $\eps^3$, because
	$\eps(1+z)\to0$ uniformly. Thus $u^{(i)}_2\ge(\xi_0/2)\eps(1+z)$. For $M_0\log\beta \le z\le M\log\beta$,  it holds that $U_{i+1}\ge\xi_0h$, and 
	$h\ge M_0 \eps\log\beta\ge(\eps+h)/2$ for large $\beta$. Since $u_{i+1}^{\ap}$ is a
	convex combination of $U_{i+1}$ and $u_2^{(i)}$, this gives (b) for $u_{i+1}^{\ap}$ near $s_i$;
	\revised{for $u_i^{\ap}$ on the left, use the lower bound for $a_1$ and the
left-hand Taylor expansion.}
For $z\le0$ one has $u_{i+1}^{\ap}=\theta_iu^{(i)}_2$.
	Since $a_2>0$, the decay of $a_2,P_2,Q_2$ gives
	\begin{equation*}
			u_{i+1}^{\ap}\ge-\eps^2\bigl(|P_2(z)|+\eps|Q_2(z)|\bigr)\ge-C\eps^2,
		\qquad |u_{i+1}^{\ap}|\le C_\tau\eps e^{-\tau|z|},
	\end{equation*}
	while (b) gives $u_i^{\ap}\ge\xi\eps(1+|z|)$ there. This proves (c) and (d)
	and (e) near $s_i$ follow for large $\beta$. For instance, for $h\le0$ one has
	$\sigma_\beta\ge\xi\eps(1+|z|)-C\eps^2\ge(\xi/2)\eps(1+|z|)$, and for $h\ge0$
	one has $|u^{(i)}_1|\le C\eps\le Cu_{i+1}^{\ap}$.
Finally, the polynomial parts of the local
	expansion differ from the segregated profiles by $O(\eps+h^2)=O(\eps)$ on
	$|h|\le M\rho$, which gives (a).
\end{proof}

Choose a sufficiently small rate
$0<\zeta\le \min\{\zeta_0, \zeta_1,\zeta_2 \}$, where $\zeta_0,  \zeta_1, \zeta_2$ are positive constants given in  Lemmas \ref{scalar-barriers-1},~\ref{two-sided-barrier},~\ref{gap-left-out-barrier} respectively.  Then
choose $\tau>\zeta$ and $0<M_0<M$ such that
\begin{equation}\label{rate-choice}
	M_0(\tau-\zeta)>1.
\end{equation}
\revised{We now estimate the error obtained by substituting $\bv u^{\ap}$
into \eqref{o-eq}.}

\begin{proposition}\label{residual}
	Let $\bv R=(R_1,R_2, \cdots,R_k)$, where 
	\begin{equation*}
		R_l=-(u_l^{\ap})''-c(u_l^{\ap})'-\mu u_l^{\ap}+\mu(u_l^{\ap})^2
	+\eps^{-3}u_l^{\ap}\sum_{j\ne l}\alpha_{lj}u_j^{\ap}.
	\end{equation*}
	All residuals vanish outside the sets $\{|r-s_i|\le M\rho\}$. On such a set,
	$R_l=0$ for $l\notin\{i,i+1\}$, and with $h=r-s_i$,
	\begin{equation}\label{Rv-active}
		|R_i| \le 	\begin{cases}
			C\eps^2(\log \beta)^2, &-M\rho\le h\le0,
			\\
			C\eps^2e^{-\zeta h/\eps},  &0\le h\le M\rho,
		\end{cases}
	\end{equation}
	\revised{For $R_{i+1}$ the polynomial bound holds for $h\ge0$, and the
	exponential bound holds for $h\le0$.} 
\end{proposition}
\begin{proof}
Outside the sets $\{|r-s_i|\le M\rho\}$ one has
	$\bv u^{\ap}=\bv U$, the products $U_lU_j$ vanish and each $U_l$ solves its
	equation, so the residuals vanish. On $\{|r-s_i|\le M\rho\}$ only $u_i^{\ap}$
	and $u_{i+1}^{\ap}$ are nonzero, so $R_l=0$ for $l\notin\{i,i+1\}$, and
	$(R_i,R_{i+1})$ is the residual of the two-component model with
	$(\alpha_{12},\alpha_{21})=(\alpha_{i,i+1},\alpha_{i+1,i})$. From now on, we
	drop the superscript $(i)$. 
 In the region $|r-s_i|\le M_0\rho$, by \eqref{a-eq},
	\eqref{P-equation} and \eqref{Q-equation}, for the first component, what remains is
	\begin{equation*}
		\begin{aligned}
				R_{1,\rm in}= & \eps^2[-cQ_1'-\mu P_1+\mu a_1^2
			+\alpha_{i,i+1}(P_1Q_2+Q_1P_2)]
			+\eps^3[-\mu Q_1+2\mu a_1P_1+\alpha_{i,i+1}Q_1Q_2]\\
			&+\mu\eps^4(P_1^2+2a_1Q_1)
			+2\mu\eps^5P_1Q_1+\mu\eps^6Q_1^2.
		\end{aligned}
	\end{equation*}
	The second component is obtained by replacing the component index by $2$ and
	$\alpha_{i,i+1}$ by $\alpha_{i+1,i}$ in the interaction terms. On
	$|z|\le M_0 \ell$ the functions $\bv a,\bv P,\bv Q$ grow at most like $|z|$,
	$z^2$, $|z|^3$, and $\eps|z|\le M_0 \eps\ell\to0$. On the side where the first
	species is present, \revised{the coefficient of $\varepsilon^2$ grows at most quadratically.
Every product of two factors with different component indices has a
rapidly decaying factor. In the remaining higher-order terms,
$\varepsilon(1+|z|)$ is bounded, so
$R_{1,\rm in}=O(\varepsilon^2(1+z^2))=O(\varepsilon^2\ell^2)$.} On the other
	side, every term contains a rapidly decaying factor, so
	$|R_{1,\rm in}|\le C_\tau\eps^2e^{-\tau|z|}$ for every fixed $\tau>0$.
	
It remains to control the interpolation regions
	$M_0\rho\le|h|\le M\rho$. Write $\Delta_1=u^{(i)}_1-U_i$, $\Delta_2=u^{(i)}_2-U_{i+1}$,
	and $\theta=\theta_i$, so that $u_i^{\ap}=U_i+\theta\Delta_1$ there. A direct computation
	gives 
	\begin{equation*}
		R_i= \theta R_{1,\rm in}-2\theta'\Delta_1'
		-(\theta''+c\theta')\Delta_1 
		-\mu\theta(1-\theta)\Delta_1^2
		-\eps^{-3}\alpha_{i,i+1}\theta(1-\theta)\Delta_1\Delta_2.
	\end{equation*}
	The identity for $R_{i+1}$ is analogous, with
	$\Delta_2$ in the derivative and square terms and $\alpha_{i+1,i}$ in the
	product term. 
	Consider the left interpolation region and put $t=|h|/\eps\in[M_0 \ell,M\ell]$. By
	Proposition \ref{matching},
	\begin{equation*}
		\begin{aligned}
				|\Delta_1|&\le C(\rho^4+\eps e^{-\tau t}),\qquad
			|\Delta_1'|\le C(\rho^3+e^{-\tau t}),\\
			|\Delta_2|&\le C\eps e^{-\tau t},\qquad
			|\Delta_2'|\le Ce^{-\tau t},
		\end{aligned}
	\end{equation*}
	while $|\theta'|\le C/\rho$ and $|\theta''|\le C/\rho^2$. For species $i$,
	which is present on the left, the polynomial parts of the terms with
	$\theta'$ and $\theta''$ are $O(\rho^2)=O(\eps^2\ell^2)$, and the square term
	is smaller. The product term satisfies
	\begin{equation}\label{product-bound}
		\eps^{-3}|\Delta_1\Delta_2|
		\le C\left(\eps^2 \ell^4 e^{-\tau t}+\eps^{-1}e^{-2\tau t}\right),
	\end{equation}
	and since $e^{-\tau t}\le\beta^{-\tau M_0}\le\eps^3$, this proves \eqref{Rv-active}. 
	For species $i+1$, which is absent on the left, the \revised{residual is estimated after multiplication by the weight} $e^{\zeta t}$. The term with
	$\theta'$ gives at most
	\begin{equation*}
		C\eps^{-1}\ell^{-1}e^{-(\tau-\zeta)t}
	\le C\eps^{-1}\ell^{-1}\beta^{-M_0(\tau-\zeta)}
	\le C\eps^2,
	\end{equation*}
	by \eqref{rate-choice}; the term with $\theta''$ is even smaller. The weighted
	right-hand side of \eqref{product-bound} is also $O(\eps^2)$ 
	because $\ell^4e^{-(\tau-\zeta)M_0\ell}\le C$, and 
	$M_0(2\tau-\zeta)>1$. The weighted square term is smaller. This proves the
	left-side bound for $R_{i+1}$, and the right interpolation region is treated in
	the same way.
\end{proof}

\subsection{Linear analysis}
\revised{
We next estimate the inverse of the linearization at $\bv u^{\ap}$.
It acts on a perturbation $\bv\phi$ by}

\begin{equation}\label{linearization}
	(\Lin\bv\phi)_j=-\phi_j''-c\phi_j'-\mu(1-2u_j^{\ap})\phi_j
	+\beta\sum_{l\ne j}\alpha_{jl}\bigl(u_l^{\ap}\phi_j+u_j^{\ap}\phi_l\bigr).
\end{equation}
\revised{We use an exponential weight where the limiting density of a component
vanishes. On its occupied interval we also weight the forcing so that its
$L^1$ norm is controlled by $\varepsilon$ times its weighted norm.} 
Fix $0<\varsigma<1$ and the rate $\zeta>0$ chosen in \eqref{rate-choice}.  For
$1\le j\le k$ define
\begin{equation}\label{weighted-norms}
	\begin{gathered}
		\varpi_{\mathcal X,j}(r)=\begin{cases}1,&r\in\bar I_j,\\
			e^{\zeta\dist(r,I_j)/\eps},&r\notin\bar I_j,\end{cases}
		\qquad
		\varpi_{\mathcal Y,j}(r)=\begin{cases}1+\bigl(\dist(r,\Sigma_j)/\eps\bigr)^{1+\varsigma},
			&r\in\bar I_j,\\
			e^{\zeta\dist(r,I_j)/\eps},&r\notin\bar I_j,\end{cases}\\
		\|\bv F\|_{\mathcal Z}=\sum_{j=1}^k\|\varpi_{\mathcal Z,j}F_j\|_\infty,
		\qquad \mathcal Z\in\{\mathcal X,\mathcal Y\}.
	\end{gathered}
\end{equation}
The norm $\nX{\cdot}$ is used for solutions $\bv\phi$ and $\nY{\cdot}$ for right-hand
sides. 
\begin{proposition}\label{inverse-linear}
There exist constants $C, \beta_0>0$, such that, for all  $\beta>\beta_0$  and every $\bv F\in C([-\uppi,\uppi])^k$, there
	is a unique Dirichlet solution of $\Lin\bv\phi=\bv F$, and
	\begin{equation}\label{inverse-estimate}
		\nX{\bv\phi}\le C\eps \nY{\bv F}, \quad \eps=\beta^{-\frac13},
	\end{equation}
\end{proposition}
\begin{proof}
	\revised{(a) We argue by contradiction. If the estimate failed, there would be}
	$\eps_\nu\to0$ and 
	\begin{equation}\label{inverse-contradiction}
		\nX{\bv\phi_\nu}=1,\qquad \gamma_\nu:=\eps_\nu\nY{\bv F_\nu}\longrightarrow0.
	\end{equation}
	We omit the index $\nu$ when no confusion arises. All components are uniformly
	bounded, and the weights give
	\begin{equation*}
		\sum_j\|F_j\|_{L^1}\le C\gamma_\nu\longrightarrow0,
	\end{equation*}
	because $\int_0^\infty(1+t^{1+\varsigma})^{-1}\dd t$ and
	$\int_0^\infty e^{-\zeta t}\dd t$ are finite. Moreover,  $|\phi_l(r)|\le e^{-\zeta\dist(r,I_l)/\eps}$ for every $l$.
	%Thus $J_i$ contains $\bar I_i\cup\bar I_{i+1}$ except
%	for a margin $\delta$ at the far ends, and in particular it contains a
%	neighbourhood of $s_i$.
	
	\revised{(b) We first pass to the limit in weighted differences near the interfaces. Let}
	$L_{\min}=\min_iL_i$, fix $0<\delta<L_{\min}/4$, and let $\nu$ be so large that
	$M\rho<\delta$. For $1\le i\le k-1$, set $J_i=[s_{i-1}+\delta,s_{i+1}-\delta]$,
	with $s_{i-1}+\delta$ replaced by $-\uppi$ if $i=1$ and $s_{i+1}-\delta$ replaced
	by $\uppi$ if $i=k-1$.  Put
	$W_i=\alpha_{i,i+1}\phi_{i+1}-\alpha_{i+1,i}\phi_i$. In the combination
	$\alpha_{i,i+1}(\Lin\bv\phi)_{i+1}-\alpha_{i+1,i}(\Lin\bv\phi)_i$, the terms
	describing the competition between species $i$ and $i+1$ cancel exactly. What remains is
	\begin{equation}\label{W-equation}
		-W_i''-cW_i'=f_i+G_i+\Gamma_i,
	\end{equation}
	where $f_i=\alpha_{i,i+1}F_{i+1}-\alpha_{i+1,i}F_i\to0$ in $L^1$,
	\begin{equation*}
			G_i=\mu\{\alpha_{i,i+1}(1-2u_{i+1}^{\ap})\phi_{i+1}
		-\alpha_{i+1,i}(1-2u_i^{\ap})\phi_i\}
	\end{equation*}
	is uniformly bounded, and $\Gamma_i$ collects the competition terms involving
	the species $l\notin\{i,i+1\}$. On $J_i$ one has $u_l^{\ap}=0$ and
	$|\phi_l|\le e^{-\zeta\delta/\eps}$ for such $l$, hence
	$|\Gamma_i|\le C\eps^{-3}e^{-\zeta\delta/\eps}\to0$ uniformly on $J_i$. \revised{The remaining forcing $f_i$ is small in $L^1$, while $G_i+\Gamma_i$
is bounded in $L^\infty$. Since $W_i$ is bounded, the mean value theorem
provides a point of $J_i$ where $|W_i'|\le C$. Integrating
$-(e^{cr}W_i')'=e^{cr}(f_i+G_i+\Gamma_i)$ from that point bounds
$W_i'$ uniformly. The integral of $f_i$ is uniformly small, and the
integral of $G_i+\Gamma_i$ is uniformly Lipschitz.} By the Arzel\`a--Ascoli theorem, after extracting a
	subsequence, $W_i\to W_i^\infty$ in $C^1(J_i)$.
	
\revised{(c) To identify the limiting profiles, use the fact that each $\phi_l$
	tends to zero away from its own interval.} on $\{\dist(r,I_l)\ge\eta\}$, for any fixed $\eta>0$,
	$\phi_l\to0$ uniformly. Hence, on $J_i\cap\{r\le s_i-\eta\}$,
	\begin{equation*}
			\phi_i=\frac{\alpha_{i,i+1}\phi_{i+1}-W_i}{\alpha_{i+1,i}}
		\longrightarrow-\frac{W_i^\infty}{\alpha_{i+1,i}},
	\end{equation*}
	and on $J_i\cap\{r\ge s_i+\eta\}$ we get
	$\phi_{i+1}\to W_i^\infty/\alpha_{i,i+1}$. Further, 
	there are functions
	$\Phi_i\in C^1(\bar I_i)$, with $\phi_i\to\Phi_i$ locally uniformly in $I_i$,
	such that $W_i^\infty=-\alpha_{i+1,i}\Phi_i$ on $J_i\cap\bar I_i$ and
	$W_i^\infty=\alpha_{i,i+1}\Phi_{i+1}$ on $J_i\cap\bar I_{i+1}$.
	We now identify the equations satisfied by the $\Phi_i$. On $J_i\cap I_i$ we
	have $u_{i+1}^{\ap}\to0$, $\phi_{i+1}\to0$ and $u_i^{\ap}\to U_i$, so passing
	to the limit in \eqref{W-equation} gives $L_{U_i}\Phi_i=0$; the same holds on
	$J_i\cap I_{i+1}$ for $\Phi_{i+1}$. Moreover $W_1(-\uppi)=0$ and
	$W_{k-1}(\uppi)=0$, so $\Phi_1(-\uppi)=\Phi_k(\uppi)=0$. Finally, $W_i^\infty$ is of
	class $C^1$ near $s_i$, which is exactly \eqref{kernel-conditions}. Lemma
	\ref{nondegeneracy} implies $\Phi_i=0$ for all $i$, hence $W_i\to0$ in $C^1(J_i)$
	for every $i$. Since the limit is the same for every subsequence, the whole
	sequence converges.
	
\revised{(d) We now use the smallness of the $W_i$ to bound the individual
components uniformly, including at the interfaces.} 
	For $2\le j\le k-1$ fixed,  choose smooth functions $0\le\eta_j^\pm\le1$ with
	$\eta_j^+=1$ on $[s_{j-1}+L_j/2,\uppi]$, $\eta_j^+=0$ on $[-\uppi,s_{j-1}+L_j/4]$,
	$\eta_j^-=1$ on $[-\uppi,s_{j-1}+L_j/2]$ and $\eta_j^-=0$ on $[s_j-L_j/4,\uppi]$. 
	 Put $\eta_1^+\equiv1$, $\eta_1^-\equiv0$, $\eta_k^+\equiv0$ and
	$\eta_k^-\equiv1$. Then $\eta_j^++\eta_j^-\ge1$, and
	$\operatorname{supp}u_j^{\ap}\cap\operatorname{supp}\eta_j^+\subset J_j$,
	$\operatorname{supp}u_j^{\ap}\cap\operatorname{supp}\eta_j^-\subset J_{j-1}$.
	Using $\alpha_{j,j+1}\phi_{j+1}=W_j+\alpha_{j+1,j}\phi_j$ where $\eta_j^+>0$,
	and $\alpha_{j,j-1}\phi_{j-1}=\alpha_{j-1,j}\phi_j-W_{j-1}$ where
	$\eta_j^->0$, the $j$th equation becomes
	\begin{equation}\label{scalar-recovery}
		[-\partial_r^2-c\partial_r+\beta S_{\beta,j}-\mu(1-2u_j^{\ap})]\phi_j
		=F_j-\beta u_j^{\ap}\bigl(\eta_j^+W_j-\eta_j^-W_{j-1}\bigr)-\beta u_j^{\ap}g_j,
	\end{equation}
	with
	\begin{equation*}
		S_{\beta,j}=\sum_{l\ne j}\alpha_{jl}u_l^{\ap}
	+\bigl(\alpha_{j+1,j}\eta_j^++\alpha_{j-1,j}\eta_j^-\bigr)u_j^{\ap}, 
	\end{equation*}
	 and
	 \begin{equation*}
	 		g_j=(1-\eta_j^+)\alpha_{j,j+1}\phi_{j+1}+(1-\eta_j^-)\alpha_{j,j-1}\phi_{j-1}
	 	+\sum_{|l-j|\ge2}\alpha_{jl}\phi_l.
	 \end{equation*}
On
	the support of $u_j^{\ap}$, we have
	$|\beta u_j^{\ap}g_j|\le C\eps^{-3}e^{-\zeta L_{\min}/(4\eps)}\to0$.
	\revised{To use scalar comparison, we need a lower bound on $S_{\beta,j}$ including near the two endpoints.} Indeed, by Lemma
	\ref{profile-geometry}, after reducing $\xi$,
	$S_{\beta,j}\ge\xi\min\{r+\uppi,\uppi-r,\eps\}$ and $|u_j^{\ap}|\le CS_{\beta,j}$. We
	can therefore apply Lemma \ref{scalar-barriers}, with
	$q_\beta=\mu(1-2u_j^{\ap})$. The second estimate in \eqref{scalar-resolvent},
	applied to the terms with $W_j$ and $W_{j-1}$, and the first one, applied to
	$F_j$ and to the terms with $g_j$, together with $\|F_j\|_\infty\le\gamma_\nu/\eps$,
	give
	\begin{equation*}
			\|\phi_j\|_\infty\le C\bigl(\|W_j\|_{L^\infty(J_j)}+\|W_{j-1}\|_{L^\infty(J_{j-1})}\bigr)
		+C\eps\gamma_\nu+C\eps^{-1}e^{-\zeta L_{\min}/(4\eps)}\longrightarrow0.
	\end{equation*}
	Write $\varkappa_\nu=\sum_j\|\phi_j\|_\infty\to0$.
	
(e) It remains to show that $\phi_j$ is small
	also when multiplied by the exponential weight outside $I_j$. Fix $j\le k-1$
	and consider the half-interval $[s_j+\eps,\uppi]$ to the right of $I_j$. There
	we keep the $j$th equation in its original form,
	\begin{equation*}
			\Big[-\partial_r^2-c\partial_r+\beta\sum_{l\ne j}\alpha_{jl}u_l^{\ap}
		-\mu(1-2u_j^{\ap})\Big]\phi_j
		=F_j-\beta u_j^{\ap}\sum_{l\ne j}\alpha_{jl}\phi_l.
	\end{equation*}
By
	\eqref{weighted-norms} and Lemma \ref{profile-geometry}(c) at a fixed rate
	$\tau>\zeta$, the right-hand side is at most
	$C\eps^{-2}(\varkappa_\nu+\eps\gamma_\nu)e^{-\zeta(r-s_j)/\eps}$. By Lemma
	\ref{profile-geometry}(e), the one-sided part of Lemma \ref{scalar-barriers-1}
	applies with $s=s_j$ and $\widetilde S_\beta=\sum_{l\ne j}\alpha_{jl}u_l^{\ap}$.
	At the left endpoint $|\phi_j|\le\varkappa_\nu$, while
	$B_\beta(s_j+\eps)\ge e^{-\zeta}/2$ for large $\nu$; at $\uppi$ both vanish.
Consequently, 
\begin{equation*}
		\sup_{s_j+\eps\le r\le\uppi}e^{\zeta(r-s_j)/\eps}|\phi_j(r)|
	\le C(\varkappa_\nu+\eps\gamma_\nu)\longrightarrow0.
\end{equation*}
	On $[s_j,s_j+\eps]$ the weight is bounded by $e^\zeta$. The region to the left
	of $I_j$ is treated by reflection, and on $\bar I_j$ the weight equals one.
	Thus $\nX{\bv\phi_\nu}\to0$, which contradicts \eqref{inverse-contradiction}
	and proves \eqref{inverse-estimate}.
	
(f) For fixed $\beta$, the Dirichlet operator is Fredholm of
	index zero: on $(H^2\cap H_0^1)^k\to(L^2)^k$ it is a compact perturbation of
	the diagonal Dirichlet Laplacian. \revised{In one dimension any kernel element is classical by the equation.
The estimate gives a trivial kernel, so the Fredholm alternative
gives a solution for each forcing. For continuous forcing the equation
then gives $C^2$ regularity.}
\end{proof}

\subsection{Exact positive solution}
We now look for an exact solution of the form $\bv u=\bv u^{\ap}+\bv\phi$,
where the correction $\bv\phi$ vanishes at $\pm\uppi$. Inserting this into
\eqref{o-eq} and using the definitions of $\bv R$ and $\Lin$ given in Proposition \ref{residual} and  \eqref{linearization} respectively, we find that
$\bv\phi$ must solve
\begin{equation}\label{perturbation}
	\Lin\bv\phi=-\bv R-\bv N(\bv\phi),
\end{equation}
and
\begin{equation*}
	N_j(\bv\phi)=\mu\phi_j^2+\beta\phi_j\sum_{l\ne j}\alpha_{jl}\phi_l,
	\qquad 1\le j\le k.
\end{equation*}
\begin{lemma}\label{nonlinear-estimate}
	For $\bv X,\bv Y\in C([-\uppi,\uppi])^k$ and large $\beta$,
	\begin{align}
		\nY{\bv N(\bv X)}&\le C\beta\nX{\bv X}^2,\label{N-bound}\\
		\nY{\bv N(\bv X)-\bv N(\bv Y)}&\le C\beta(\nX{\bv X}+\nX{\bv Y})\nX{\bv X-\bv Y}.\label{N-lipschitz}
	\end{align}
\end{lemma}
\begin{proof}
	Let $A=\nX{\bv X}$, so that $|X_l(r)|\le Ae^{-\zeta\dist(r,I_l)/\eps}$ for
	every $l$. Consider $X_jX_l$, $l\ne j$, measured with the weight
	$\varpi_{\mathcal Y,j}$. If $r\in\bar I_j$, then  $\dist(r,I_l)\ge\dist(r,\Sigma_j)$, and the
	product times the weight is bounded by $CA^2$ because
	$\sup_{t\ge0}(1+t^{1+\varsigma})e^{-\zeta t}<\infty$. If $r\notin\bar I_j$,
	the factor $|X_j|$ compensates the weight. For the square terms
	$\mu X_j^2$, outside $\bar I_j$ the weight is compensated by one of the two
	factors, while on $\bar I_j$ the weight is at most
	$C\eps^{-(1+\varsigma)}\le C\beta$. This proves \eqref{N-bound}. Writing the
	differences of squares and of products as sums of terms of the same type
	proves \eqref{N-lipschitz}.
\end{proof}

\begin{proposition}\label{existence}
	For some fixed $\Lambda>0$ and all sufficiently large $\beta$,
	\eqref{perturbation} has a unique solution in
	\begin{equation*}
		\mathcal B_\beta=\{\bv X\in C([-\uppi,\uppi])^k:
		\bv X(\pm\uppi)=0,\quad \nX{\bv X}\le\Lambda\eps^3(\log\beta)^{3+\varsigma}\}.
	\end{equation*}
	The solution is smooth and, in particular,
	\begin{equation}\label{correction-bound}
		\nX{\bv\phi}\le \Lambda \beta^{-1}(\log\beta)^{3+\varsigma}.
	\end{equation}
\end{proposition}
\begin{proof}
	Set $\varsigma_*=3+\varsigma$, $\ell=\log\beta$ and
	$T(\bv X)=-\Lin^{-1}(\bv R+\bv N(\bv X))$, so that the solutions of
	\eqref{perturbation} are the fixed points of $T$. By Propositions
	\ref{residual}, \ref{inverse-linear} and Lemma \ref{nonlinear-estimate},
	\begin{equation*}
			\nX{T(\bv X)}\le C\eps\bigl(\eps^2\ell^{\varsigma_*}+\beta\nX{\bv X}^2\bigr)
		=C\eps^3\ell^{\varsigma_*}+C\eps^{-2}\nX{\bv X}^2.
	\end{equation*}
 For $\bv X \in 	\mathcal B_\beta$, the second
	term is at most $C\Lambda^2\eps^4\ell^{2\varsigma_*}$, which is much smaller
	than  $\eps^3\ell^{\varsigma_*}$  because $\eps\ell^{\varsigma_*}\to0$. Hence, choosing 
	$\Lambda$ large, $T$
	maps the ball $\mathcal B_\beta$ into itself for large $\beta$. By \eqref{N-lipschitz}, its
	Lipschitz constant on the ball is at most $C\Lambda\eps\ell^{\varsigma_*}<1$
	for large $\beta$. Banach's fixed point theorem gives a unique fixed point in
	the ball. The equation and one-dimensional regularity give smoothness.
\end{proof}

\begin{proposition}\label{positivity}
	The solution given by Proposition \ref{existence} has all components
	positive throughout $(-\uppi,\uppi)$ for all sufficiently large $\beta$.
\end{proposition}
\begin{proof}
	Let $d_\beta=\Lambda\eps^3(\log\beta)^{3+\varsigma}=o(\eps)$, then
	$|\phi_l(r)|\le d_\beta e^{-\zeta\dist(r,I_l)/\eps}$ for every $l$.
We claim that, for every $j$,
	\begin{equation}\label{positive-active}
		u_j(r)\ge\xi_1\min\{\eps+\dist(r,\Sigma_j),\,r+\uppi,\,\uppi-r\},\qquad r\in\bar I_j,
	\end{equation}
	for a fixed $\xi_1>0$. On the region where the minimum is at least $\eps$, \eqref{positive-active} follows from
	Lemma \ref{profile-geometry}(b) directly, since $|\phi_j|\le d_\beta=o(\eps)$. Otherwise
	$r$ is close to $-\uppi$ (and $j=1$) or to $\uppi$ (and $j=k$). Consider the
	latter case; the former is symmetric. 
Choose a fixed $d>0$ with $3d<L_k$. On
	$[\uppi-2d,\uppi]$ the approximate solution equals $\bv U$, so $u_l=\phi_l$ is
	exponentially small for $l\ne k$, and
	\begin{equation*}
			-\phi_k''-c\phi_k'=\mu(1-2U_k)\phi_k-\mu\phi_k^2
		-\beta u_k\sum_{l\ne k}\alpha_{kl}u_l,
		\qquad \|\phi_k\|_\infty+\|\phi_k''+c\phi_k'\|_\infty\le C d_\beta.
	\end{equation*}
Consequently 
	$\|\phi_k'\|_{L^\infty(\uppi-d,\uppi)}\le Cd_\beta$. Since $\phi_k(\uppi)=0$, it
	follows that $|\phi_k(r)|\le Cd_\beta(\uppi-r)$, while Hopf's lemma gives
	$U_k\ge c_1(\uppi-r)$ near $\uppi$. This proves \eqref{positive-active}.
	
 Fix $1\le j\le k-1$ and consider
	$[s_j,\uppi]$, to the right of $I_j$. There $u_j$ solves the linear equation
	\begin{equation*}
			\Bigl[-\partial_r^2-c\partial_r+\beta\sum_{l\ne j}\alpha_{jl}u_l
		-\mu(1-u_j)\Bigr]u_j=0,\qquad u_j(s_j)\ge\xi_1\eps,\quad u_j(\uppi)=0.
	\end{equation*}
 We claim that  on $[s_j,\uppi]$
	\begin{equation}\label{potential-lower}
		\beta\sum_{l\ne j}\alpha_{jl}u_l\ge\beta\xi_2\min\{\eps,\uppi-r\}-Q_1
	\end{equation}
	with fixed $\xi_2,Q_1>0$. Every $r\in[s_j,\uppi]$ lies in some $\bar I_m$ with
	$m\ne j$, and by \eqref{positive-active} the term with $l=m$ is at least
	$\alpha_{jm}\xi_1\min\{\eps,\uppi-r\}$. Other component $u_l$ can be negative at $r$ only if
	$r\notin\bar I_l$. If $\dist(r,I_l)\le M\rho$, then $r$ lies within $M\rho$ of a
	interface point of $I_l$, thus
	$u_m\ge\xi_1\eps$, while $u_l\ge-C\eps^2-d_\beta$ by Lemma
	\ref{profile-geometry}(c); the negative part is therefore much smaller than
	the positive one. If $\dist(r,I_l)\ge M\rho$, then $u_l=\phi_l\ge-d_\beta
	\beta^{-\zeta M} $, and $ \beta d_\beta\beta^{-\zeta M}=\Lambda (\log\beta)^{3+\varsigma} 	\beta^{-\zeta M} $ is bounded. This proves
	\eqref{potential-lower}.
	Furthermore,
the barrier $1-e^{-(\uppi-r)/\eps}$ gives the
	maximum principle on $[s_j,\uppi]$,  as in Lemma
	\ref{scalar-barriers-1} with $\zeta=0$ (applied with $s=s_j-\eps$); the bounded
	quantities $Q_1$ and $\mu(1-u_j)$ are included in $q_\beta$. Thus $u_j\ge0$
	on $[s_j,\uppi]$. An interior zero would be a minimum, with zero derivative,
	and uniqueness for the linear initial value problem would force $u_j\equiv0$,
	a contradiction. Hence $u_j>0$ there. By reflection, $u_j>0$ on
	$[-\uppi,s_{j-1}]$ for $j\ge2$, and together with \eqref{positive-active} all
	components are positive in $(-\uppi,\uppi)$.
\end{proof}

\subsection{\revised{Proof of Theorem \ref{main-theorem}}}
\begin{proof}[Completion of the proof of Theorem \ref{main-theorem}]
	 Proposition \ref{existence} gives a classical solution of \eqref{o-eq}, and Proposition
	\ref{positivity} shows that all its components are positive. Lemma
	\ref{profile-geometry}(a) and \eqref{correction-bound} prove
	the upper bound in \eqref{global-convergence}. 	At any fixed interface point $s_i$, 
	\begin{equation*}
			u_{\beta,i}(s_i)=\eps a_1^{(i)}(0)+O(\eps^2)+O(d_\beta).
	\end{equation*}
	Since $a_1^{(i)}(0)>0$, $d_\beta=o(\eps^2)$ and $U_i(s_i)=0$, this
	gives the lower bound in \eqref{global-convergence}.
	
	To prove the estimate \eqref{main-inner-expansion}, we fix an interface point $s_i$ and drop the superscript
	$(i)$. In the region $|h|\le M_0\rho$ the approximate solution equals the local
	expansion \eqref{inner}. There $|\bv P(z)|\le C(1+z^2)$ and
	$|\bv Q(z)|\le C(1+|z|^3)$, hence
	\begin{equation*}
			|\eps^2\bv P(z)+\eps^3\bv Q(z)|\le C(\eps^2+h^2).
	\end{equation*}
	In the interpolation regions, the same estimate holds by \eqref{matching-final}.
	Since the correction satisfies $d_\beta=o(\eps^2)$, this proves
	\eqref{main-inner-expansion} for $|h|\le M\rho$. For $M\rho\le|h|\le M'\rho$ one
	has $(u_i,u_{i+1})=(U_i,U_{i+1})+O(d_\beta)$. On the side where species $i$ is
	present, $U_i-\eps a_1(z)=(U_i+p_ih)-\eps(a_1+p_iz)=O(h^2)+O(\eps^3)$ by
	\eqref{Vtaylor} and Lemma \ref{decay}, while on the other side
	$\eps a_1(z)=O(\eps^3)$; the same holds for $u_{i+1}$. For $j\notin\{i,i+1\}$
	and $|h|\le M'\rho$ one has $u_j^{\ap}=0$ and $\dist(r,I_j)\ge L_{\min}/2$ for
	large $\beta$, so $|u_j|=|\phi_j|\le d_\beta e^{-\zeta L_{\min}/(2\eps)}$.

	Finally, on each fixed compact set of values of $z$,  the cutoff equals one for
	large $\beta$, so the rescaled solution $\widehat u_l(z)=\eps^{-1}u_l(s_i+\eps z)$
	converges uniformly to $\bv a$ for $l\in\{i,i+1\}$ and to zero for the other
	species. In the rescaled variables, the equations read
	\begin{equation*}
			-\widehat u_i''-c\eps\widehat u_i'
		-\mu\eps^2\widehat u_i+\mu\eps^3\widehat u_i^2
		+\widehat u_i\sum_{l\ne i}\alpha_{il}\widehat u_l=0,
	\end{equation*}
	and similarly for $\widehat u_{i+1}$. Their coefficients are locally
	bounded, so local estimates for ordinary differential equations upgrade the
	uniform convergence to convergence in $C^2$ on compact sets, which proves
	\eqref{main-inner-limit}.
\end{proof}

\section{\revised{The problem with repeated species labels}} \label{1-d-multiplicity}
The
interval $(-\uppi,\uppi)$ is now divided into $m> k$ subintervals, and a
labeling tells us which species occupies each of them. 
When $m=k$ and $\sigma$ is the identity, this is Definition~\ref{def-profile}.
The flux balance \eqref{flux-sigma} is the same as \eqref{flux}, with the
coefficients of the two species that actually meet at $s_l$.   The number of labelings of $m$ intervals by $k$ species is equal to 
\begin{equation*}
	\sum_{i=0}^{k}(-1)^i \binom ki(k-i)(k-i-1)^{m-1}
\end{equation*}
by the inclusion-exclusion principle.

\begin{proposition}\label{sigma-existence}
	Let $\sigma$ be a labeling of $m$ intervals. A segregated profile defined in Definition \ref{def-labelled-profile} with
	labeling $\sigma$ exists if and only if $\mu-c^2/4>m^2/4$, and it is then
	unique. Moreover, it is nondegenerate in the sense of Lemma
	\ref{nondegeneracy}, with $\alpha_{i,i+1}$ and $\alpha_{i+1,i}$ replaced by
	$\alpha_{\sigma(l),\sigma(l+1)}$ and $\alpha_{\sigma(l+1),\sigma(l)}$.
\end{proposition}
\begin{proof}
	The proof of Proposition \ref{outer-existence} uses only the lengths of the
	intervals and the ratios in the flux balance. With
	$\vartheta_l=\alpha_{\sigma(l+1),\sigma(l)}/\alpha_{\sigma(l),\sigma(l+1)}$ in
	place of $\vartheta_i$, the recursion \eqref{length-recursion} and the
	shooting argument apply without change.
	Likewise, the proof of Lemma \ref{nondegeneracy} uses only the positivity of the coefficients appearing in the
	conditions at the interface points; it does not use the fact that different species occupy different
	intervals.
\end{proof}
Let $\sigma$ be a fixed labelling of $m$ intervals with
\revised{$\mu-c^2/4>m^2/4$, and $\bv U$ the corresponding segregated profile.
If species $j$ occupies several intervals, a component $(a,b)$ of
the complement of $\Omega_j$ may have interfaces at both ends.
Such a gap is a positive distance from the physical endpoints.
The linear and positivity arguments must therefore allow data at both
ends of the gap. We use Lemma~\ref{two-sided-barrier} for this purpose.} 

\subsection{The approximate solution}
\revised{At each interface point $s_l$}, we use the model profile and the corrections of
Section~\ref{sec-inner} with
$(\alpha_{12},\alpha_{21},p,q)=(\alpha_{\sigma(l),\sigma(l+1)},
\alpha_{\sigma(l+1),\sigma(l)},p_l,q_l)$, and we define the local expansion
$\bv u^{(l)}$ by \eqref{inner}; its first component approximates species
$\sigma(l)$ and its second species $\sigma(l+1)$. Proposition \ref{matching}
holds without change, with $V_l$ and $V_{l+1}$ in place of $U_i$ and
$U_{i+1}$. \revised{Use $\theta_l(r)=\chi((r-s_l)/\rho)$ to join the local expansions
to the logistic solutions. The contribution from $I_l$ is}
\begin{equation*}
	V_l^{\ap}=(1-\theta_{l-1}-\theta_l)V_l+\theta_{l-1}u^{(l-1)}_2+\theta_l  u^{(l)}_1,
\end{equation*}
with $V_l$ extended by zero, and the $j$th component of the approximate
solution is
\begin{equation*}
	u_j^{\ap}=\sum_{\sigma(l)=j} V_l^{\ap}.
\end{equation*}
Since two intervals of the same species are never adjacent, the functions
$V_l^{\ap}$ with $\sigma(l)=j$ have disjoint supports. Two species are present
at the same point only near an interface point where they meet.
Lemma
\ref{profile-geometry} holds with $I_l$ replaced by $\Omega_l$ in (b), (c)
and (e). 
Proposition \ref{residual} holds with the same proof. Indeed, near $s_l$ only
the species $\sigma(l)$ and $\sigma(l+1)$ are present, and the computation is
that of the two-component model; away from the interface points the
approximate solution is the segregated profile.  
\subsection{The linear analysis}
\revised{In the weights \eqref{weighted-norms}, replace $\bar I_j$ by
$\Omega_j$ and use the interfaces $\Sigma_j$ of species $j$.
We prove the estimate of Proposition~\ref{inverse-linear} by the
same contradiction argument, including the case where a species occurs
on both sides of another species.}

(a) For each interface point $s_l$,  set
$W_l=\alpha_{\sigma(l),\sigma(l+1)}\phi_{\sigma(l+1)}
-\alpha_{\sigma(l+1),\sigma(l)}\phi_{\sigma(l)}$, and let $J_l$ be defined as
before, with $s_{l-1}$ and $s_{l+1}$. On $J_l$ only the species $\sigma(l)$
and $\sigma(l+1)$ are present in the approximate solution. Every other species
$j$ occupies only intervals that are disjoint from $I_l\cup I_{l+1}$, so its
distance from $J_l$ is at least $\delta$, and $\phi_j$ is exponentially small
there. Hence \eqref{W-equation} holds on $J_l$, and $W_l$ converges in
$C^1(J_l)$ along a subsequence. 

(b) We also note that
$\phi_{\sigma(l+1)}$ is exponentially small on $J_l\cap\{r\le s_l-\eta\}$ for
every fixed $\eta>0$.  Symmetrically, $\phi_{\sigma(l)}$ is exponentially
small on $J_l\cap\{r\ge s_l+\eta\}$. On each interval $I_l$, the component $\phi_{\sigma(l)}$
converges locally uniformly to a function $\Phi_l\in C^1(\bar I_l)$, which is
obtained from $W_{l-1}$ and $W_l$ exactly as before. The functions
$\Phi_1,\dots,\Phi_m$ are indexed by the intervals, not by the species, and
they satisfy the conditions of Lemma \ref{nondegeneracy} with the coefficients
of Proposition \ref{sigma-existence}. Hence they vanish, and $W_l\to0$ in
$C^1(J_l)$ for every $l$.

(c) To prove $\phi_j$ tends to zero uniformly, we work separately on each interval
$I_l\subset\Omega_j$. One new situation must be taken into account: the two
neighbors of $I_l$ may be occupied by the same species. We therefore choose, for each $l$, smooth
functions $0\le\eta_l^\pm\le1$ with \revised{$\eta_l^+ + \eta_l^-=1$ at distance at most $\delta$ from $\bar I_l$}, $\eta_l^+=0$ on $[-\uppi,s_{l-1}+L_l/4]$ and
$\eta_l^-=0$ on $[s_l-L_l/4,\uppi]$; for $l=1$ we take $\eta_1^+=1$ and
$\eta_1^-=0$, and for $l=m$ we take $\eta_m^+=0$ and \revised{$\eta_m^-=1$. On the}
support of $V_l^{\ap}$ we write each coupling term as
$\alpha_{ji}\phi_i=\eta_l^+\alpha_{ji}\phi_i+\eta_l^-\alpha_{ji}\phi_i$. If
$i=\sigma(l+1)$, we rewrite the first part by means of $W_l$; if
$i=\sigma(l-1)$, we rewrite the second part by means of $W_{l-1}$. \revised{On the support of $V_l^{\ap}$, each remaining term involves a
component $\phi_i$ evaluated at distance at least $L_{\min}/4$ from
its own occupied set $\Omega_i$, by the choice of cutoffs. Hence
these terms are exponentially small.} Summing over the intervals of
species $j$, we obtain an equation of the form \eqref{scalar-recovery}, in
\revised{which the coefficient of $V_l^{\ap}$ in the new potential is}
$\alpha_{\sigma(l+1),j}\eta_l^++\alpha_{\sigma(l-1),j}\eta_l^-\ge
\min_{i\ne j}\alpha_{ij}$ near $I_l$. Lemma \ref{scalar-barriers} then shows
that $\phi_j\to0$ uniformly, as before.

(d) The complement of $\Omega_j$ consists of at most one interval
to the left of $\Omega_j$, at most one interval to the right, and the gaps.
On the first two we use the one-sided barriers of Lemma
\ref{scalar-barriers-1}, exactly as before. On each gap $(a,b)$ we use Lemma
\ref{two-sided-barrier} on $[a+\eps,b-\eps]$, with
$\widetilde S_\beta=\sum_{l\ne j}\alpha_{jl}u_l^{\ap}$. The assumption
$\eps^2\beta\widetilde S_\beta\ge\xi$ holds by Lemma
\ref{profile-geometry}(e), because the gap is at a positive distance from
$\pm\uppi$. As in the proof of Proposition \ref{inverse-linear}, the
right-hand side is bounded by
$C\eps^{-2}(\varkappa_\nu+\eps\gamma_\nu)e^{-\zeta\dist(r,\Omega_j)/\eps}$, and at the
endpoints of $[a+\eps,b-\eps]$ we have $|\phi_j|\le\varkappa_\nu$, while
$B_\beta^{\mathrm{gap}}\ge e^{-\zeta}$. Since $e^{-\zeta\dist(r,\Omega_j)/\eps}\le B_\beta^{\mathrm{gap}}$, comparison with
$C(\varkappa_\nu+\eps\gamma_\nu)B_\beta^{\mathrm{gap}}$ gives
\begin{equation*}
	\sup_{a+\eps\le r\le b-\eps}e^{\zeta\dist(r,\Omega_j)/\eps}|\phi_j(r)|
	\le C(\varkappa_\nu+\eps\gamma_\nu)\longrightarrow0.
\end{equation*}
On $[a,a+\eps]$ and $[b-\eps,b]$ the weight is bounded by $e^\zeta$. This
completes the contradiction argument.

\subsection{The nonlinear problem and positivity}
Lemma \ref{nonlinear-estimate} holds with the same proof: if $r\in\Omega_j$,
then any path from $r$ to $\Omega_l$, $l\ne j$, crosses a point of $\Sigma_j$,
so $\dist(r,\Omega_l)\ge\dist(r,\Sigma_j)$. Proposition \ref{existence} then
follows with the same proof. In the proof of Proposition \ref{positivity}, the positivity
on $\Omega_j$ is obtained interval by interval, as before. Outside $\Omega_j$,
the regions to the left and to the right of $\Omega_j$ are treated as before,
and on each gap $(a,b)$ the function $u_j$ solves a linear equation whose
potential is at least $\beta\xi_2\eps-Q_1$ on $[a,b]$, by the argument leading
to \eqref{potential-lower}, with positive values at $a$ and $b$. Since this
potential is of order $\eps^{-2}$ on the whole closed gap, the constant
function $1$ is a positive strict supersolution there, and the maximum
principle gives $u_j\ge0$ on the gap; an interior zero is excluded as in
the proof of Proposition \ref{positivity}, so $u_j>0$ there. \revised{At each interface, the approximate solution has exactly two nonzero
components. The distance from that interface to the closed occupied set of
any other species is bounded below independently of $\beta$. Thus the
estimates at the end of the proof of Theorem~\ref{main-theorem} give
the same rate and local expansion here. Since the set of admissible
labellings for $k\le m\le m_0$ is finite and their profiles are distinct,
their minimum pairwise uniform distance is positive whenever there are
at least two profiles. Taking $\beta$ sufficiently large gives distinct
solutions for all these profiles at once. This proves
Theorem~\ref{main-theorem-sigma}.}

\section{\revised{The radial problem}} \label{radial-sys}

Write
\begin{equation*}
	\D v=v''+\frac{n-1}{r}v'=r^{1-n}(r^{n-1}v')',\qquad
	\varrho(r)=r^{n-1},\qquad \varepsilon=\beta^{-1/3}.
\end{equation*}
The problem \eqref{o-system-radial} is 
\begin{equation*}
	\begin{cases}
		-\D u_i=\mu u_i(1-u_i)-\beta u_i\displaystyle\sum_{j\ne i}\alpha_{ij}u_j,
		&0<r<R,\\
		u_i'(0)=0,\quad u_i(R)=0,&i=1,2,\cdots,k,
	\end{cases}
\end{equation*}
where $\mu>0$ and $\alpha_{ij}>0$ for $i\ne j$.  A classical radial solution means that $x\mapsto u_i(|x|)$ is a classical solution on $B_R\subset\R^n$. Let $\upnu=(n-2)/2$ and let $j_{\upnu,l}$ be the $l$th positive zero of $J_\upnu$. By Lemma \ref{eigenvalue--radial}, the radial eigenvalues of the Dirichlet Laplacian on $B_R$ are
\begin{equation*}
	\lambda_l^\rad(B_R)=\frac{j_{\upnu,l}^2}{R^2},\qquad l=1,2,\ldots,
\end{equation*}
with eigenfunctions $\phi_l(x)=\phi_l(r)= \revised{r^{-\upnu}J_\upnu}(  \frac{j_{\upnu,l}}{R} r)$, $r=|x|$.
\revised{For each labelling, the $m$-th radial eigenvalue is the threshold for
the existence of the corresponding segregated profile.}
\begin{proposition}\label{profile-radial}
	For every fixed labelling $\sigma$, a segregated profile in Definition~\ref{def-profile-radial} exists if and only if
	\begin{equation*}
		 \mu>\lambda_m^\rad(B_R)=j_{(n-2)/2,m}^2/R^2.
	\end{equation*}
	When it exists, it is unique, and its radial linearization has a trivial kernel. 
\end{proposition}
\begin{proof}
	\revised{Fix the weights $d_l$ of \eqref{eq-d-l} and let $b_m(t)$ be the
$m$-th positive zero of the solution started at $Y(0)=t$.
Lemmas~\ref{return} and \ref{range} show that $b_m$ is strictly
increasing and has range $(j_{\upnu,m}/\sqrt\mu,\infty)$.
Thus $b_m(t)=R$ has exactly one solution in the stated range of $\mu$.
On each shell set $V_l=(-1)^{l-1}Y/d_l$. These functions solve the
logistic equation, and continuity of $Y'$ is precisely the flux
balance \eqref{flux-condi} with the weights \eqref{eq-d-l}.
Conversely, every radial segregated profile gives a solution of
this initial-value problem with $t=V_1(0)\in(0,1)$, which proves
uniqueness.} 
		For completeness, necessity also follows directly from the spectrum. For a profile, $\lambda_m(-\D-\mu+\mu W)=0$, while strict potential monotonicity gives
		\begin{equation*}
			0>\lambda_m(-\D-\mu)=\lambda_m^\rad(B_R)-\mu.
		\end{equation*}
		The nondegeneracy assertion is Corollary~\ref{transmission}.
	\end{proof}

\subsection{Model profile near an interface point}
\revised{Fix an interface $s=s_l>0$. With $h=r-s$ and $z=h/\varepsilon$,
we approximate the densities of species $\sigma(l)$ and $\sigma(l+1)$ by}
\begin{equation}\label{inner-profile-radial}
	\bx v^\mathrm{in}(r)=\eps\bx a(z)+\eps^2\bx P(z)+\eps^3\bx Q(z), \quad \eps=\beta^{-\frac13}.
\end{equation}
The leading terms of adjacent species behave like the profile $\bv a$ in Proposition \ref{inner-profile} with $(\alpha_{12},\alpha_{21},p,q)$ replaced by $(\alpha_{\sigma(l),\sigma(l+1)},\alpha_{\sigma(l+1),\sigma(l)},p_l,q_l)$.  The linearization at $\bv a$ is given by \eqref{inner-operator} with $\bv b$ replaced by $(\alpha_{\sigma(l),\sigma(l+1)},\alpha_{\sigma(l+1),\sigma(l)})^T$. % Lemma \ref{inverse} still holds. 
 Furthermore, the corrections differ because of the presence of $(n-1)/r$ in the operator $-\D$. 
  Near a fixed positive $s$,
 \begin{equation*}
 	\frac{n-1}{s+\eps z}=c_s+\eps c_s'z+O(\eps^2z^2).
 \end{equation*}
  Let $s\in(0, R)$ be a fixed interface. Set
\begin{equation*}
	c_s=\frac{n-1}{s},\qquad c_s'=-\frac{n-1}{s^2},\qquad
	\gamma_s=\mu+c_s'-c_s^2=\mu-\frac{n(n-1)}{s^2}.
\end{equation*}
\revised{The constant part $c_s$ contributes as in the interval problem.
The variation of the radial drift contributes at the next order.
Write $\bv b=(\alpha_{\sigma(l),\sigma(l+1)},\alpha_{\sigma(l+1),\sigma(l)})^T$.}
\begin{proposition}\label{radial-correctors}
	Near an interface point $s=s_l$, 
	define
	\begin{equation*}
		\bv P(z)=-\frac{c_s}{10}\bigl(z^2\bv a'(z)+4z\bv a(z)\bigr).
	\end{equation*}
	Then
	\begin{equation*}
		\begin{aligned}
			\begin{cases}
				\M\bv P =c_s \bv a', \quad  
				\Dif{\bv P}=\frac{c_s p}{2 b_1 }z^2,
				\\
				(P_1, P_2)= (c_spz^2/2 +\mathcal E_-, \mathcal E_-) \  \mbox{at } -\infty,
				\quad
				(P_1, P_2)=(	\mathcal E_+, -c_sqz^2/2+\mathcal E_+)
				\ \mbox{at } +\infty.
			\end{cases}
		\end{aligned}
	\end{equation*}
	Moreover, there is a unique smooth $\bv Q$  such that
	\begin{equation*}
		\begin{aligned}
			\begin{cases}
				\M\bv Q=c_s\bv P'+c_s'z \bv a'+\mu\bv a-\bv b\,P_1P_2, \quad
				\Dif{\bv Q}=\frac{\gamma_s p}{6b_1 }z^3,
				\\
				(Q_1, Q_2)= \big( 	\gamma_s pz^3/6 +\mathcal E_- , \mathcal E_-\big)
				\  \mbox{at } -\infty, \quad 	(Q_1, Q_2)= \big(
				\mathcal E_+ ,  - \gamma_s qz^3/6+\mathcal E_+
				\big) \ \mbox{at } +\infty.
			\end{cases}
		\end{aligned}
	\end{equation*}
%	The corrections are independent of $\beta$.  Here each occurrence of $\mathcal E_\pm$ denotes a possibly different
	%function that is rapidly decaying at the corresponding end, together with
%	its derivatives.
\end{proposition}

\begin{proof}
\revised{Insert \eqref{inner-profile-radial} into the two equations and expand
$(n-1)/(s+\varepsilon z)$ at $s$. The terms of order
$\varepsilon^{-1}$ give \eqref{a-eq}. Those of order one give
$\M\bv P=c_s\bv a'$, so the computation in
Proposition~\ref{correctors} applies with $c_s$ in place of $c$.
At order $\varepsilon$ the drift contributes $c_s\bv P'+c_s'z\bv a'$.
The logistic term contributes $\mu\bv a$, and the product of the first
corrections contributes $-\bv bP_1P_2$ on the right-hand side.
If $\bv H$ denotes this right-hand side, then
\[
\Dif{\bv H}
=c_s\Dif{\bv P}'+c_s'z\Dif{\bv a}'+\mu\Dif{\bv a}
=-\gamma_s pz/b_1 .
\]
The required polynomial difference therefore satisfies
\eqref{D-compatible}. By the decay of $\bv a$ and $\bv P$,
$H_2$ decays rapidly at the left end and $H_1$ at the right end.
Lemma~\ref{inverse} gives $\bv Q$ and its stated asymptotics.}
\end{proof}

\revised{The same coefficients appear in the Taylor expansions of the
segregated densities.} At the interface point, the left profile has $V(s)=0$, $V'(s)=-p$, and
\begin{equation*}
V''(s)=c_sp,\qquad V'''(s)=\gamma_sp.
\end{equation*}
The right profile has the opposite slope sign. Therefore
\begin{equation*}
	\begin{aligned}
			V_\ell(s+h)&=-ph+\frac{c_sp}{2}h^2+\frac{\gamma_sp}{6}h^3+O(h^4),\quad h\le0,
			\\
		V_{\ell+1}(s+h)&=qh-\frac{c_sq}{2}h^2-\frac{\gamma_sq}{6}h^3+O(h^4), \quad h\ge0.
	\end{aligned}
\end{equation*}
\revised{These Taylor polynomials agree with the polynomial parts of
\eqref{inner-profile-radial} at the respective ends, by
Propositions~\ref{inner-profile} and \ref{radial-correctors}.} 

\subsection{Approximate solution} 
\revised{Near $s_l$, use the cutoff $\theta_l(r)=\chi((r-s_l)/\rho)$ to interpolate
between the two logistic densities and their local expansion. The two
components of the approximate solution there are}
\begin{equation*}
	\revised{(1-\theta_l)(V_l,V_{l+1})+
	\theta_l\bx v^{\mathrm{in},l}},
\end{equation*}
where the outer functions are extended by zero to the other side.  Every other component is zero there. Away from interface neighborhoods, use $\bx U$. This defines $\bx u^\ap$. The supports of two occurrences of the same species are disjoint. The result is smooth, satisfies the radial endpoint conditions, and equals $\bx U$ near $0$ and $R$. 
\revised{The Taylor coefficients and the polynomials in the local expansion
agree through degree three. Subtracting these polynomials leaves the
Taylor remainder for the logistic density and an exponentially decaying
term from the expansion, just as in Proposition~\ref{matching}.}
\begin{lemma}
	On an interpolation region $M_0 \rho \le|h|\le M\rho$, the difference $\Delta$ between \revised{the local expansion for the species occupying that side and its logistic profile} satisfies, for $j=0,1,2$,
	\begin{equation*}
		|\partial_r^j\Delta|\le C|h|^{4-j}+C_\tau\eps^{1-j}e^{-\tau|h|/\eps}.
	\end{equation*}
	\revised{For the other component the bound is $C_\tau\eps^{1-j}e^{-\tau|h|/\eps}$.} The estimate holds for every fixed $\tau>0$.
\end{lemma}

\revised{The lower bounds follow as in Lemma~\ref{profile-geometry}. In each
interface neighborhood the leading profile of the species present on a
given side is bounded below by a positive multiple of $\varepsilon+|h|$.
The corrections are smaller there because $\varepsilon(1+|z|)$ tends
uniformly to zero. Away from the interfaces the approximation equals the
segregated profile. Its first component is positive at the center and its
last component has a strict Hopf slope at $R$.}
\begin{lemma}
	There exist constants $\xi,C>0$ such that for every sufficiently large $\beta$, 
	\begin{equation*}
		\begin{aligned}
			&	\|\bx u^\ap-\bx U\|_\infty  \le C\eps,\qquad \|\bx u^\ap\|_\infty\le C,
			\\
		&	u_l^\ap(r) \ge \xi \min\{\eps+\dist(r,\Sigma_l),R-r\},
			\quad r\in\Omega_l,
			\\
	&		|u_l^\ap(r)| \le C_\tau\eps e^{-\tau\dist(r,\Omega_l)/\eps},\quad  r\notin\Omega_l,
	\\
	&		u_l^\ap \ge-C\eps^2,\qquad
			\sum_l u_l^\ap\ge \xi \min\{\eps,R-r\}.
		\end{aligned}
	\end{equation*}
	\revised{The component $u_l^\ap$ is zero when $\dist(r,\Omega_l)\ge M\rho$.} Outside $\Omega_l$,
	\begin{equation*}
		\revised{\sum_{j\ne l}\alpha_{lj}u_j^\ap}\ge \xi \min\{\eps,R-r\},\qquad
		|u_l^\ap|\le C\sum_{j\ne l} u_j^\ap.
	\end{equation*}
%	All statements hold for repeated labels.
\end{lemma}

Define
\begin{equation*}
	R_i=-\D u_i^\ap-\mu u_i^\ap+\mu(u_i^\ap)^2
	+\beta u_i^\ap\sum_{j\ne i}\alpha_{ij} u_j^\ap.
\end{equation*}
We use the norms \eqref{weighted-norms} with $\bar I_j$ replaced by $\Omega_j$
and $\dist(r,I_j)$ replaced by $\dist(r,\Omega_j)$. Consequently, \revised{the radial drift gives two kinds of additional terms. First,
$(n-1)/r-c_s-c_s'h=O(h^2)$ multiplies the radial derivative of the
local expansion. Second, $c_s'h$ multiplies the derivatives
$\varepsilon\bv P'$ and $\varepsilon^2\bv Q'$, which contribute at
order $\varepsilon^2$ or higher. On the side where a species is
present, the derivatives of $\bv a,\bv P,\bv Q$ grow at most like
$1$, $1+|z|$, $1+z^2$, respectively. Since $\varepsilon|z|$ is small,
both kinds of terms are $O(\varepsilon^2(1+z^2))$. On the other side
all these derivatives decay rapidly, giving the exponential estimate.
The cutoff terms are bounded as in Proposition~\ref{residual}
because the radial drift is bounded near each fixed $s_l>0$.}
\begin{proposition}
	\revised{The residual vanishes outside the interface neighborhoods and satisfies}
	\begin{equation*}
		\nY{\bx R}\le C\eps^2(\log \beta)^{3+\varsigma}.
	\end{equation*}
\end{proposition}

\subsection{Linear analysis}
The linearization at the approximation is
\begin{equation*}
	(\Lr_\eps\bx\phi)_i=-\D\phi_i-\mu(1-2u_i^\ap)\phi_i
	+\eps^{-3}\sum_{j\ne i}\alpha_{ij}(u_j^\ap\phi_i+u_i^\ap\phi_j),
\end{equation*}
with $\phi_i'(0)=0$, $\phi_i(R)=0$.

\begin{proposition}\label{inverse-radial}
	For sufficiently small $\eps$ and every continuous radial forcing $\bx F$, there is a unique regular radial solution of $\Lr_\eps\bx\phi=\bx F$, and
	\begin{equation*}
		\nX{\bx\phi}\le C\eps\nY{\bx F}.
	\end{equation*}
\end{proposition}
\begin{proof}
\revised{Suppose the estimate fails. As in
\eqref{inverse-contradiction}, take a sequence with
$\nX{\bv\phi}=1$ and $\gamma_\varepsilon=\varepsilon\nY{\bv F}\to0$.
The weights give $\sum_j\|F_j\|_{L^1(0,R)}\le C\gamma_\varepsilon$.
For each interface use $W_l$ as in
Section~\ref{1-d-multiplicity}, with $J_1$ starting at $0$ and $J_{m-1}$
ending at $R$. The mutual competition cancels and
$-\D W_l=f_l+G_l+\Gamma_l$ on $J_l$, where $f_l\to0$ in $L^1$,
$G_l$ is bounded and $\Gamma_l\to0$ uniformly. The last assertion follows
from the positive distance of $J_l$ to the occupied sets of species
other than the two at $s_l$.

On each $J_l$ away from the center, integration after multiplication
by $r^{n-1}$ bounds $W_l'$ and gives relative compactness in $C^1$,
as in the proof of Proposition~\ref{inverse-linear}. At the center we
use instead
\[
W_1'(r)=-\int_0^r(t/r)^{n-1}(f_1+G_1+\Gamma_1)(t)\,\dd t.
\]
The contribution of $f_1$ is bounded by $\|f_1\|_{L^1}$.
The contribution of the bounded terms is $O(r)$ and is uniformly
Lipschitz, as follows by differentiating the integral for $r>0$.
Consequently, $W_1$ is relatively compact in $C^1(J_1)$ and its
limit has zero derivative at $0$.

On each region where a species is positive the other component in the local difference
vanishes in the limit away from the interface. Passing to the limit in
the difference equations thus gives $L_{V_l}\Phi_l=0$ on each shell.
The $C^1$ limits give the transmission conditions
\eqref{transmission-r}, and the limits at the center and at $R$ give
$\Phi_1'(0)=0$ and $\Phi_m(R)=0$. Corollary~\ref{transmission} forces
every $\Phi_l$ to vanish, hence $W_l\to0$ uniformly on each $J_l$.

To estimate each component, split the coupling terms using the cutoffs
in part (c) of Section~\ref{1-d-multiplicity}. On the support of the
approximation from the $l$th shell, substitute the differences
$W_{l-1}$ and $W_l$ in the corresponding terms. The terms left over
involve components exponentially small on that support. The scalar
operator for species $j$ has a potential
$\varepsilon^{-3}S_{\varepsilon,j}-q_{\varepsilon,j}$, where
$q_{\varepsilon,j}=\mu(1-2u_j^{\ap})$ is bounded. On each shell of
species $j$, the coefficient of its approximate density in
$S_{\varepsilon,j}$ is bounded below by a positive constant.
The lower bounds for the approximate densities then give
\[
S_{\varepsilon,j}\ge \xi\min\{\varepsilon,R-r\},
\qquad |u_j^{\ap}|\le C S_{\varepsilon,j}.
\]
Indeed, on a region where a species is positive its own density supplies the lower
bound. Elsewhere a different species supplies it, and any negative
approximate density is $O(\varepsilon^2)$ and is supported near an
interface, where a present species has density at least
$\xi\varepsilon$. Applying the second estimate in \eqref{scalar-inv}
to the terms with $W_{l-1}$ and $W_l$ bounds their contribution by
$C\max_l\|W_l\|_{L^\infty(J_l)}$. The first estimate bounds the
contributions of $F_j$ and the exponentially small terms by $o(1)$
since $\|F_j\|_\infty\le\gamma_\varepsilon/\varepsilon$. Thus
$\varkappa_\varepsilon=\sum_j\|\phi_j\|_\infty\to0$.

Outside $\Omega_j$ we return to the original $j$th equation. With
$d(r)=\dist(r,\Omega_j)$, its right-hand side is bounded by
$C\varepsilon^{-2}(\varkappa_\varepsilon+
\varepsilon\gamma_\varepsilon)e^{-\zeta d(r)/\varepsilon}$, using the
exponential estimate on $u_j^{\ap}$. The scalar potential satisfies the
lower bounds in Lemma~\ref{gap-left-out-barrier}. Use its barrier near
the center, on each bounded gap, or near $R$, according to the component
of the complement. On a gap the sum of the two exponentials is
comparable to $e^{-\zeta d(r)/\varepsilon}$. The center barrier has
the same property and zero derivative at $0$. The outer barrier
vanishes at $R$ and its image under the operator controls the forcing
up to $R$. At points a distance $\varepsilon$ from an interface,
these barriers are bounded below by a positive constant, so they
also dominate the boundary data after multiplication by
$C(\varkappa_\varepsilon+\varepsilon\gamma_\varepsilon)$.
Within distance $\varepsilon$ of an interface the weight itself is
bounded. We conclude that $\nX{\bv\phi}\to0$, a contradiction.

For fixed $\varepsilon$ the radial operator on the ball is a compact
perturbation of the diagonal Dirichlet Laplacian, acting from the radial
subspace of $(H^2\cap H^1_0)^k$ to radial $(L^2)^k$. The estimate
gives a trivial kernel, and the Fredholm alternative gives existence.
Elliptic regularity gives continuity of the solution for continuous
forcing, and the radial integral formula gives $C^2$ regularity in $r$,
including at $0$.}
\end{proof}

Write $\bx u=\bx u^\ap+\bx\phi$. The correction equation is
\begin{equation}\label{fixed-point-eq}
	\Lr_\eps\bx\phi=-\bx R-\bx N(\bx\phi),\qquad
	N_i(\bx\phi)=\mu\phi_i^2+\eps^{-3}\phi_i\sum_{j\ne i}\alpha_{ij}\phi_j.
\end{equation}

\begin{proposition}\label{contraction-radial}
	\revised{There is a fixed $\Lambda>0$} such that, for all sufficiently large $\beta$, \eqref{fixed-point-eq} has exactly one solution in
	\begin{equation*}
		\mathcal B_\eps=\{\bx\phi\in C([0,R])^k:\bx\phi(R)=0,
		\ \nX{\bx\phi}\le\Lambda\eps^3(\log \beta)^{3+\varsigma}\}.
	\end{equation*}
\revised{The correction is regular at the origin and satisfies $\phi_i'(0)=0$.
	The solution $\bv u^{\ap}+\bv\phi$ is smooth on the ball, and the
	correction obeys}
	\begin{equation*}
		\nX{\bx\phi}\le \Lambda \eps^3(\log \beta)^{3+\varsigma}=o(\eps^2).
	\end{equation*}
\end{proposition}

\begin{proof}
\revised{The quadratic estimates in Lemma~\ref{nonlinear-estimate}
hold with $\Omega_j$ in place of $\bar I_j$. In fact, if
$r\in\Omega_j$ and $i\ne j$, a segment from $r$ to $\Omega_i$
crosses an interface of species $j$, so
$\dist(r,\Omega_i)\ge\dist(r,\Sigma_j)$.
The residual and inverse estimates therefore bound the fixed-point map by
$C\varepsilon^3(\log\beta)^{3+\varsigma}
+C\varepsilon^{-2}\nX{\bv\phi}^2$.
For a sufficiently large fixed $\Lambda$ this maps the stated ball
into itself. Its Lipschitz constant on the ball is
$O(\varepsilon(\log\beta)^{3+\varsigma})$, which tends to zero.
The contraction theorem gives the unique correction there.
The resulting solution of the nonlinear system is smooth on the ball
by elliptic regularity.}
\end{proof}

\begin{proposition}\label{positivity-radial}
	For sufficiently large $\beta$, the solution in
	Proposition~\ref{contraction-radial} has every component strictly
	positive on $[0,R)$.
\end{proposition}
\begin{proof}
\revised{Write $d_\beta=\Lambda\varepsilon^3(\log\beta)^{3+\varsigma}$.
On a region where a species is positive, the lower bound for $u_j^{\ap}$ dominates the
correction $|\phi_j|\le d_\beta=o(\varepsilon)$ except possibly
within distance $\varepsilon$ of $R$. On a fixed neighborhood of $R$,
the other species are exponentially small and the approximation
equals the segregated profile. The correction equation on that
neighborhood gives $|\D\phi_{\sigma(m)}|\le Cd_\beta$.
Using the equation and the boundary value at $R$ as in the proof of
Proposition~\ref{positivity}, we obtain
$|\phi_{\sigma(m)}(r)|\le Cd_\beta(R-r)$ near $R$.
The strict Hopf slope of $V_m$ gives the required lower bound there.

Outside the occupied set of species $j$, a species that occupies the
point contributes at least $\xi\min\{\varepsilon,R-r\}$ to
$\sum_{i\ne j}\alpha_{ji}u_i$. Near an interface, a negative component
is bounded in magnitude by $C\varepsilon^2+d_\beta$, and the positive
contribution is at least $\xi\varepsilon$. Away from the cutoff
neighborhoods of its own occupied set, a negative component has
magnitude at most $d_\beta e^{-\zeta M\rho/\varepsilon}$.
After multiplication by $\beta$, the latter bound stays bounded.
It follows that the scalar equation for $u_j$ off $\Omega_j$ has
potential at least $\xi_1\varepsilon^{-3}
\min\{\varepsilon,R-r\}-C$.

On a bounded gap or a region adjoining the center this potential is
strictly positive for small $\varepsilon$. The constant function $1$
is a strict supersolution, with zero derivative at the center.
On a region adjoining $R$, the outer barrier in
Lemma~\ref{gap-left-out-barrier} with $\zeta=0$ gives the maximum
principle up to $R$. The same differentiation applies on the full
interval from its fixed left endpoint to $R$. The interface values are
positive, and the exterior boundary condition is either zero value at
$R$ or zero derivative at the center. Scalar comparison thus gives
$u_j\ge0$ off $\Omega_j$. The strong maximum principle, applied on the
ball for the region containing the center, gives strict positivity in
the interior since the solution has positive interface values.}
\end{proof}

\begin{proof}[Completion of the proof of Theorem~\ref{thm-main-radial}]
\revised{The approximate solution is within $C\varepsilon$ of
$\bv U$, and the correction is $o(\varepsilon^2)$ in the uniform norm.
Evaluating at any interface gives
$u_{\beta,\sigma(l)}(s_l)=\varepsilon a_1^{(l)}(0)+O(\varepsilon^2)$
with $a_1^{(l)}(0)>0$, which proves the matching lower bound.

Where the cutoff is one, the bounds
$|\bv P(z)|\le C(1+z^2)$ and $|\bv Q(z)|\le C(1+|z|^3)$
give the stated local error. In the interpolation region the Taylor
and matching estimates give the same bound. Beyond that region, up to
$|h|\le M'\rho$, the segregated density differs from
$\varepsilon\bv a(h/\varepsilon)$ by $O(h^2)$ plus a rapidly
decaying term. The distance to the occupied sets of the remaining
species is bounded below throughout this neighborhood, so their
weighted correction estimate gives the claimed exponential bound.
On compact $z$-intervals the rescaled drift
$\varepsilon(n-1)/(s_l+\varepsilon z)$ tends uniformly to zero.
The rescaled equations and the uniform convergence of the densities
give convergence in $C^2$ on compact sets.

For each $m$ with $k\le m\le M_\mu$, Proposition~\ref{profile-radial}
gives one profile for each admissible labelling. Inclusion--exclusion
counts them as in Section~\ref{1-d-multiplicity}. Since there are only
finitely many such profiles, uniform convergence to their distinct
limits gives the stated lower bound for all sufficiently large
$\beta$.}
\end{proof}

\section{\revised{A vanishing competition coefficient}} \label{vanish-rate}
\begin{proof}[Proof of Proposition \ref{prop-zero}]
\revised{
Let $I$ be an interval or shell occupied by species $l$. In the
one-dimensional problem write $A=-\partial_r^2-c\partial_r$ and
$\omega=e^{cr}$. In the radial problem write $A=-\D$ and
$\omega=r^{n-1}$. The logistic equation on $I$ implies
$\lambda_1(I)<\mu$. By shrinking the Dirichlet endpoints slightly,
choose $D\Subset D'\Subset I$ with $\lambda_1(D)<\mu$.
If $I$ contains the center, keep the regular condition there and
shrink only the outer radius. The strict inequality persists because
the first eigenvalue is continuous under these changes of endpoints.

Uniform convergence gives $u_l\ge c_0>0$ on $D'$ for large $\beta$.
For $h\notin\{j,l\}$, positivity and the signs of the interaction
coefficients imply
\[
(A+\beta\alpha_{hl}c_0-\mu)u_h\le0\quad\text{on }D'.
\]
The densities are uniformly bounded. On an interval or an annulus
with endpoints $a,b$, compare with a constant times
$e^{-\kappa(r-a)}+e^{-\kappa(b-r)}$, with $\kappa$ a sufficiently
small multiple of $\sqrt\beta$. The zeroth-order term dominates the
second derivatives and the bounded drift, so this is a supersolution.
On a ball of radius $b$, use
$\cosh(\kappa r)/\cosh(\kappa b)$. Its derivative is zero at the
center, and its radial Laplacian is bounded above by $n\kappa^2$
times the function. In either case comparison gives
$u_h\le C e^{-c_*\sqrt\beta}$ on $D$. In particular
$\beta u_h\to0$ uniformly there.

Let $\varphi>0$ be the first eigenfunction on $D$. Use
$\alpha_{jl}=0$, multiply the equation for $u_j$ by
$\omega\varphi$ and integrate by parts. Since
$A=-\omega^{-1}(\omega\partial_r)'$, this gives
\begin{equation}\label{zero-contr}
\int_D\omega
\big[\mu-\lambda_1(D)-\mu u_j
-\beta\sum_{h\notin\{j,l\}}\alpha_{jh}u_h\big]
u_j\varphi\,\dd r
=[\omega u_j\varphi']_{\partial D}\le0.
\end{equation}
The boundary term is nonpositive by the Hopf signs at the
Dirichlet endpoints. It vanishes at the center by regularity.
Because $U_j=0$ on $D$, we have $u_j\to0$ there. The decay
just proved for the other species makes the bracket strictly
positive for large $\beta$, contradicting \eqref{zero-contr}.}

\end{proof}

\appendix
\section{\revised{Scalar results on an interval}}\label{1d-shooting-maximum}
\subsection{\revised{The logistic equation and its endpoint slopes}}
 \begin{lemma}\label{eigenvalue-n}
	The eigenvalue  problem 
	\begin{equation*}
		\begin{cases}
			-y''-cy'=\lambda y & x\in (a,b)
			\\
			y(a)=0=y(b)
		\end{cases}
	\end{equation*}
	admits eigenvalues $\lambda_j=\big(\frac{j \uppi}{b-a}\big)^2+\frac{c^2}{4}$  with eigenfunctions $\phi_j(x)=\mathrm{e}^{-\frac{c}{2}x}\sin \frac{j\uppi(x-a)}{b-a}$,  where $ j=1,2,\cdots$. 
\end{lemma}
\begin{proof}
	See \cite[p. 86]{1990-ode}.  Indeed, the substitution
	$w=e^{-cx/2}y$ gives $-w''-cw'=e^{-cx/2}(-y''+c^2y/4)$. 
\end{proof}

\begin{proposition}\label{exist-unique}
	Equation \eqref{1d-1} admits a unique positive solution if and only if
	$\mu>\lambda_1(a,b)$. This solution is smooth up to the
	endpoints, and satisfies $0<w<1$, $w'(a)>0$ and $w'(b)<0$.
\end{proposition}
\begin{proof}
	
\revised{
Let $\phi_1>0$ be the first Dirichlet eigenfunction. For any positive
solution, multiplication by $e^{cx}\phi_1$ and integration by parts
give
\[
(\mu-\lambda_1)\int_a^b e^{cx}w\phi_1\,\dd x
=\mu\int_a^b e^{cx}w^2\phi_1\,\dd x>0.
\]
Thus $\mu>\lambda_1$. Conversely, for this range a sufficiently small
multiple of $\phi_1$ is a subsolution, and the constant $1$ is a
supersolution. The sub- and supersolution method gives a positive
solution between them
\cite[Theorem~9.4]{amann-1976}. If it reached $1$ at an interior
point, its derivative there would be zero. Uniqueness for the
initial-value problem would then give the constant solution $1$,
in conflict with the boundary values. The equation gives smoothness,
and Hopf's lemma gives the endpoint signs.

To prove uniqueness, let $v$ and $w$ be positive solutions. Their
nonzero Hopf slopes imply that $v^2/w$ has zero trace and is an
admissible test function. Subtracting the equation tested against
$v^2/w$ for $w$ from the equation tested against $v$ for $v$, with
weight $e^{cx}$ in both integrals, gives
\[
\int_a^b e^{cx}\left(v'-\frac{v w'}{w}\right)^2\,\dd x
=\mu\int_a^b e^{cx}(w-v)v^2\,\dd x.
\]
Add the identity with $v$ and $w$ interchanged. The left-hand side
is nonnegative and the right-hand side equals
$-\mu\int_a^b e^{cx}(v-w)^2(v+w)\,\dd x$. Hence $v=w$.}

\end{proof}

\begin{lemma}\label{scalar-stability}
	For a positive solution $w$ of \eqref{1d-1}, the linearized Dirichlet operator
	\begin{equation*}
			L_w=-\frac{\mathrm{d}^2}{\mathrm{d}x^2} - c \frac{\mathrm{d}}{\mathrm{d}x} -\mu(1-2w)
	\end{equation*}
	has a strictly positive first eigenvalue. In particular, it is invertible
	and satisfies the maximum principle.
\end{lemma}
\begin{proof}
	Using the equation of $w$, one finds $L_w w=\mu w^2>0$. Let $\lambda$ and
	$\varphi>0$ be the principal eigenvalue and eigenfunction of $L_w$. Testing
	$L_w\varphi=\lambda\varphi$ by $e^{cx}w$ and using the self-adjointness
	of $L_w$ in the weighted space, we obtain
	\begin{equation*}
			\lambda\int_a^b e^{cx}w\varphi\dd x
		=\int_a^b e^{cx}\varphi\,L_ww\dd x
		=\mu\int_a^b e^{cx}w^2\varphi\dd x>0,
	\end{equation*}
	so $\lambda>0$. The maximum principle follows from \cite[Theorem B.4]{struwe-2008}. Fredholm theory yields the invertibility.
\end{proof}

\begin{lemma}\label{slope-monotonicity}
	Assume $\mu>c^2/4$, and set $L_* =\uppi/\sqrt{\mu-c^2/4}$.
	For the positive solution $w_L$ on $(0,L)$, $L>L_*$, define
	\begin{equation*}
			B(L)=w_L'(0),\qquad A(L)=-w_L'(L).
	\end{equation*}
	These functions are smooth, $A'(L),B'(L)>0$, and
	$A(L),B(L)\to0$ as $L\downarrow L_*$.
\end{lemma}
\begin{proof}
	Rescaling $(0, L)$ to a fixed interval, Lemma \ref{scalar-stability} and the
	implicit function theorem show that $w_L$ depends smoothly on $L$. Let
	$z=\partial_L w_L$ be the derivative with respect to the length. Then
	\begin{equation*}
			L_{w_L}z=0,\qquad z(0)=0,\quad z(L)=-w_L'(L)>0.
	\end{equation*}
	By the maximum principle in Lemma \ref{scalar-stability}, $z>0$ in $(0,L)$,
	and Hopf's lemma gives $z'(0)>0$. Since $B'(L)=z'(0)$, this proves
	$B'(L)>0$. The reflection $x\mapsto L-x$ exchanges the two endpoints and
	replaces $c$ by $-c$, so the same argument gives $A'(L)>0$.
	
	Furthermore, as $L\downarrow L_*$, the bounds $0<w_L<1$ and the equation give
	compactness in $C^2$ after rescaling to a fixed interval. Any limit is a
	nonnegative solution on an interval of length $L_*$. A nonzero limit would
	be positive by the strong maximum principle, contradicting
	Proposition~\ref{exist-unique}, since $\mu=\lambda_1$ at the threshold. Thus
	the limit is zero, together with its endpoint derivatives.
\end{proof}

For $1\le i\le k$,  let $E_i^-,E_i^+\in C^2(\bar I_i)$ solve
\begin{equation}\label{Epm}
	L_{U_i}E_i^\pm=0\ \text{in }I_i,\qquad
	E_i^-(s_{i-1})=1,\ E_i^-(s_i)=0,\qquad
	E_i^+(s_{i-1})=0,\ E_i^+(s_i)=1,
\end{equation}
where $U_i$ is the $i$-th component of the segregated profile in Definition \ref{def-profile}.  
These problems are uniquely solvable by Lemma \ref{scalar-stability}. Both
functions are positive in $I_i$ by the maximum principle, and Hopf's lemma
gives
\begin{equation}\label{Epm-hopf}
	(E_i^-)'(s_i)<0<(E_i^+)'(s_{i-1}).
\end{equation}
The next lemma compares the slopes of $E_i^\pm$ at the other endpoint. It
gives the sign of the diagonal entries of the matrix that appears in the proof of Lemma~\ref{nondegeneracy}.

\begin{lemma}\label{transversality}
	For every $i$, we have
	\begin{equation*}%\label{transversality-ineq}
		(E_i^+)'(s_i)<-c<(E_i^-)'(s_{i-1}).
	\end{equation*}
\end{lemma}
\begin{proof}
	Notice that  $L_{U_i}U_i'=0$,  moreover, at the zeros of $U_i$,  the
	equation gives $U_i''=-cU_i'$.
	Consider $E_i^--U_i'/q_{i-1}$. It solves the linear equation, vanishes at
	$s_{i-1}$ and equals $p_i/q_{i-1}>0$ at $s_i$. By 
	Lemma \ref{scalar-stability}, it is positive in $I_i$, and Hopf's lemma at
	$s_{i-1}$ gives
	$(E_i^-)'(s_{i-1})-U_i''(s_{i-1})/q_{i-1}>0$. Since
	$U_i''(s_{i-1})=-cq_{i-1}$, this is $(E_i^-)'(s_{i-1})+c>0$.
	Similarly, $E_i^++U_i'/p_i$ vanishes at $s_i$ and is positive at $s_{i-1}$,
	hence it is positive in $I_i$ and has negative derivative at $s_i$. Since
	$U_i''(s_i)=cp_i$, this is $(E_i^+)'(s_i)+c<0$.
\end{proof}

\subsection{Maximum principles for linear analysis and positivity}
\begin{lemma}\label{scalar-barriers}
	Let $S_\beta\ge0$ be continuous on $[-\uppi,\uppi]$ with
$S_\beta(r)\ge\xi\min\{r+\uppi,\uppi-r,\eps\}$ with $\xi>0$ fixed, and  $q_\beta$ be any real continuous function with
	$\|q_\beta\|_\infty\le\qbar$, \revised{uniformly in $\beta$.}
 Then, for sufficiently large $\beta$,
	the Dirichlet operator
	\begin{equation*}
		A_\beta=-\partial_r^2-c\partial_r+\beta S_\beta-q_\beta
	\end{equation*}
	has the maximum principle and is invertible, with
	\begin{equation}\label{scalar-resolvent}
		\|A_\beta^{-1}f\|_\infty\le C\eps^2\|f\|_\infty,\qquad
		\|A_\beta^{-1}(\beta S_\beta g)\|_\infty\le C\|g\|_\infty.
	\end{equation}
\end{lemma}
\begin{proof}
(a) Set
\begin{equation*}
		b_\beta(r)=\bigl(1-e^{-(r+\uppi)/\eps}\bigr)
	\bigl(1-e^{-(\uppi-r)/\eps}\bigr),\qquad 0<b_\beta\le1.
\end{equation*}
	Write 
	$t=\min\{r+\uppi,\uppi-r\}/\eps$. Direct calculations yield
	\begin{equation*}
			\eps^2\big[ -\partial_r^2-c\partial_r+\beta S_\beta \big]b_\beta
		\ge (1-|c|\eps)(e^{-(r+\uppi)/\eps}+e^{-(\uppi-r)/\eps})+\xi\min\{t,1\}b_\beta\ge b_*>0.
	\end{equation*}
	Indeed, if $t\le1$ the first term has a positive lower bound, and if $t\ge1$,
	then $b_\beta\ge(1-e^{-1})^2$ and the second term does. For large $\beta$, it follows that  $A_\beta b_\beta\ge(b_*/2)\eps^{-2}$.
	Thus 
	$\tilde b_\beta: =b_\beta+\delta_0$, with $\delta_0=\min\{1,b_*/(4\qbar)\}$
	(any $\delta_0\in(0,1]$ if $\qbar=0$), satisfies $\tilde b_\beta\ge\delta_0$ on
	$[-\uppi,\uppi]$ and $A_\beta\tilde b_\beta\ge(b_*/2)\eps^{-2}-\delta_0\qbar>0$.
	
	Now suppose $A_\beta w\ge0$ and $w(\pm\uppi)\ge0$ but $\min w<0$, and let
	$t_0=\max(-w/\tilde b_\beta)>0$. Then $w+t_0\tilde b_\beta\ge0$ vanishes at an
	interior point, where its first derivative vanishes and its second
	derivative is nonnegative, so $A_\beta(w+t_0\tilde b_\beta)\le0$ there. This
	contradicts $A_\beta(w+t_0\tilde b_\beta)\ge t_0A_\beta\tilde b_\beta>0$.  Therefore, the maximum principle holds. 
	In
	particular, the homogeneous Dirichlet problem has only the trivial solution,
	and the two-point boundary value problem is uniquely solvable.
	
	(b) Comparing $\pm A_\beta^{-1}f$ with
	$C\eps^2\|f\|_\infty b_\beta$ gives the first estimate in
	\eqref{scalar-resolvent}. For the second, choose a fixed $C_0$ so large that
	\begin{equation*}
			A_\beta(1+C_0\eps^2b_\beta)
		=\beta S_\beta-q_\beta+C_0\eps^2 A_\beta b_\beta\ge\beta S_\beta.
	\end{equation*}
	This function is positive at both endpoints, so comparing
	$\pm A_\beta^{-1}(\beta S_\beta g)$ with $\|g\|_\infty(1+C_0\eps^2b_\beta)$ gives
	the claim.
\end{proof}

\begin{lemma}\label{scalar-barriers-1}
		 Let 
	\revised{$J_\beta=[s+\eps,\uppi]$ with $s\in(-\uppi,\uppi)$.
Assume $q_\beta$ is continuous with $|q_\beta|\le\qbar$ and
$\widetilde S_\beta$ is continuous on $J_\beta$ with} $\eps^2\beta\widetilde S_\beta(r)\ge\xi\min\{1,(\uppi-r)/\eps\}$.
	Then, there exist $\zeta_0>0$ and $b_0>0$, depending only on $\xi$, $c$
	and $\qbar$, such that, for $0\le\zeta\le\zeta_0$, the function
	\begin{equation*}
		B_\beta(r)=e^{-\zeta(r-s)/\eps}\bigl(1-e^{-(\uppi-r)/\eps}\bigr)
	\end{equation*}
	satisfies
	\begin{equation}\label{weighted-barrier-ineq}
		\bigl[-\partial_r^2-c\partial_r+\beta\widetilde S_\beta-q_\beta\bigr]B_\beta
		\ge b_0\eps^{-2}e^{-\zeta(r-s)/\eps}\quad\text{on }J_\beta.
	\end{equation}
	The operator in brackets has the Dirichlet maximum principle on
	$J_\beta$. Reflected statements hold near the left endpoint.
\end{lemma}
\begin{proof}
	Set $\hat s=(r-s)/\eps$, $t=(\uppi-r)/\eps$ and
	$b(t)=1-e^{-t}$. A direct differentiation gives
	\begin{align*}
		\eps^2e^{\zeta\hat s}
		[-\partial_r^2-c\partial_r+\beta\widetilde S_\beta-q_\beta]B_\beta
		& \ge (1-2\zeta-|c|\eps)e^{-t}
		+[\xi\min\{1,t\}-\zeta^2-|c|\eps\zeta-\qbar\eps^2](1-e^{-t}).
	\end{align*}
	The function $e^{-t}+\xi\min\{1,t\}(1-e^{-t})$ has a positive infimum on
	$[0,\infty)$. Choosing $\zeta_0$ sufficiently small and then $\beta$
	sufficiently large proves \eqref{weighted-barrier-ineq}. At $\zeta=0$, the
	function $B_\beta+\delta_0$, with $\delta_0>0$ small, is bounded below by
	$\delta_0$ and is mapped to a function bounded below by
	$b_0\eps^{-2}-\delta_0\qbar>0$, so the argument in Lemma \ref{scalar-barriers} gives the Dirichlet
	maximum principle on $J_\beta$. Reflection only changes $c$ to $-c$.
\end{proof}

\begin{lemma}\label{two-sided-barrier}
	Let $G_\beta=[a+\eps,b-\eps]$, $-\uppi<a<b<\uppi$ with $b-a\ge4\eps$, 
	and let $\widetilde S_\beta$ and $q_\beta$ be
	continuous on $G_\beta$ with $\eps^2\beta\widetilde S_\beta\ge\xi$ and
	$|q_\beta|\le\qbar$, where $\xi>0$ and $\qbar\ge0$ are fixed. Then there exists $\zeta_1>0$, depending only on $\xi$ and $c$,
	such that, for $0\le\zeta\le\zeta_1$ and $\beta$ large the function
	\begin{equation*}
			B_\beta^{\mathrm{gap}}(r)=e^{-\zeta(r-a)/\eps}+e^{-\zeta(b-r)/\eps}
	\end{equation*}
	satisfies
	\begin{equation*}
		\bigl[-\partial_r^2-c\partial_r+\beta\widetilde S_\beta-q_\beta\bigr]
	B_\beta^{\mathrm{gap}}\ge\frac\xi2\eps^{-2}B_\beta^{\mathrm{gap}}
	\quad\text{on }G_\beta.
	\end{equation*}
	In particular, the operator in brackets has the Dirichlet maximum principle on
	$G_\beta$.
\end{lemma}
\begin{proof}
	For $f(r)=e^{-\zeta(r-a)/\eps}$ we have $f'=-(\zeta/\eps)f$ and
	$f''=(\zeta/\eps)^2f$, hence
	\begin{equation*}
			\bigl[-\partial_r^2-c\partial_r+\beta\widetilde S_\beta-q_\beta\bigr]f
		=\Bigl(-\frac{\zeta^2}{\eps^2}+\frac{c\zeta}{\eps}+\beta\widetilde S_\beta
		-q_\beta\Bigr)f
		\ge\eps^{-2}\bigl(\xi-\zeta^2-|c|\zeta\eps-\qbar\eps^2\bigr)f.
	\end{equation*}
	The same bound holds for $e^{-\zeta(b-r)/\eps}$, with $c$ replaced by $-c$ in
	the middle term. Choosing $\zeta_1$ with $\zeta_1^2+|c|\zeta_1\le\xi/4$ and then
	$\beta$ large gives the inequality. Since $B_\beta^{\mathrm{gap}}$ is positive on
	the closed interval $G_\beta$, the maximum principle follows exactly as in the
	proof of Lemma \ref{scalar-barriers}.
\end{proof}

\section{\revised{Scalar radial results}}\label{radial-shooting-maximum}
\subsection{Properties of the operator and shooting argument}

 \begin{lemma}\label{eigenvalue--radial}
	The eigenvalue  problem 
	\begin{equation*}
		\begin{cases}
			-v''-\frac{n-1}{r}v'=\lambda v, & r\in (0,R)
			\\
			v'(0)=0=v(R)
		\end{cases}
	\end{equation*}
	admits eigenvalues $\lambda_l= ( \frac{j_{\upnu,l}}{R}  )^2 $  with eigenfunctions $\phi_l(r)=r^{-\upnu}J_\upnu(  \frac{j_{\upnu,l}}{R} r)$,  where  $\upnu =\frac{n-2}{2}$, $ l=1,2,\cdots$. 
\end{lemma}
\begin{proof}
 Indeed, the substitution
	$v(r)= r^{-\frac{n-2}{2}} w(\sqrt{\lambda} r)$ gives $x^2 w''+x w' +( x^2- \upnu^2)w=0$ with $w(\sqrt{\lambda} R)=0$.
	\revised{Regularity at $0$ selects $w(x)=J_{\upnu}(x)$ up to a constant}, $\sqrt{\lambda}_l R =j_{\upnu,l}$ and then $\phi_l(r)=r^{-\upnu}J_\upnu(  \frac{j_{\upnu,l}}{R} r)$. Recall the series definition of $J_{\upnu}$ in \cite[Chapter III]{Watson-1952}
	\begin{equation*}
		J_{\upnu}(z)= \sum_{m=0}^{\infty} \frac{(-1)^m}{m! \Gamma(m+\upnu+1)} \big( \frac{z}{2} \big)^{2m+\upnu},
	\end{equation*}
	we obtain
	\begin{equation*}
		\begin{aligned}
					\phi_l(r)=&r^{-\upnu}J_\upnu(  \frac{j_{\upnu,l}}{R} r) =
					\sum_{m=0}^{\infty} \frac{(-1)^m}{m! \Gamma(m+\upnu+1)} \big( \frac{j_{\upnu,l}}{2R} \big)^{2m+\upnu} r^{2m}, 
			\\
				\phi_l'(r) =& 
					\sum_{m=1}^{\infty} \frac{(-1)^m (2m)}{m! \Gamma(m+\upnu+1)} \big( \frac{j_{\upnu,l}}{2R} \big)^{2m+\upnu} r^{2m-1}.
		\end{aligned}
	\end{equation*}
	Accordingly, $	\phi_l'(0)=0$ and $	\phi_l(R)=0$.
\end{proof}

For a regular solution of $-\D v=f$, integration gives
\begin{equation*}
	v'(r)=-r^{1-n}\int_0^r t^{n-1}f(t)\dd t.
\end{equation*}
In particular $v'(0)=0$ and $v''(0)=-f(0)/n$ when $f$ is continuous. 
The radial  space on $B_b$ is
\begin{equation*}
	H_\rad^1(0,b)=\left\{v:\int_0^b\varrho(v'^2+v^2)\dd r<\infty,\quad v(b)=0\right\}, \quad  \varrho=r^{n-1}
\end{equation*}
understood as the radial subspace of $H_0^1(B_b)$. For a bounded continuous radial potential $q$, the regular radial realization of $- \D + q$ is self-adjoint, has a compact resolvent and simple eigenvalues, and its $j$-th eigenfunction has exactly $j-1$ zeros in $(0,b)$. 
On an annulus $(a,b)$, $a>0$, use zero traces at both endpoints.

\begin{lemma}\label{scalar-radial}
	On a ball $(0,b)$ with the regular condition at $0$, or on an annulus $(a,b)$ with two Dirichlet conditions, the problem
	\begin{equation*}
			-\D V=\mu V(1-V)
	\end{equation*}
	has a positive solution if and only if $\mu>\lambda_1$, the first eigenvalue with the stated endpoint conditions. It is unique, satisfies $0<V<1$, and has strict Hopf slopes at all Dirichlet endpoints. Its linearized operator
	\begin{equation*}
		L_V=-\D-\mu+2\mu V
	\end{equation*}
	has a positive first eigenvalue. The positive solution depends smoothly on admissible endpoints, with the origin fixed for balls.
\end{lemma}
\begin{proof}
	\revised{(a) Let $\varphi_1>0$ be the first eigenfunction. Multiply the
equation by $\varrho\varphi_1$ and integrate by parts, using the
regularity condition at the center on a ball.}
	\begin{equation*}
		(\mu-\lambda_1)\int\varrho V\varphi_1\dd r
	=\mu\int\varrho V^2\varphi_1\dd r>0.
	\end{equation*}
	Conversely, a sufficiently small multiple of $\varphi_1$ is a subsolution and the constant $1$ is a supersolution. The sub- and supersolution method on the ball or annulus gives a positive solution. The maximum principle gives $V\le1$, and equality at an interior maximum would force the solution to be identically $1$, contradicting the Dirichlet boundary. The Hopf lemma supplies the endpoint signs.
	
	\revised{(b) For uniqueness, let $V,W>0$ be two solutions. The Hopf slopes
make $V^2/W$ an admissible test function at the Dirichlet endpoints,
and it is regular at the center on a ball. Test the equation of $V$
by $V$, and that of $W$ by $V^2/W$, using weight $\varrho$. Subtract
the two identities to get}
	\begin{equation*}
			\int\varrho\left(V'-\frac{VW'}W\right)^2\dd r
		=\mu\int\varrho(W-V)V^2\dd r.
	\end{equation*}
	Add the identity with $V, W$ interchanged. The right side becomes
	$-\mu\int\varrho(V-W)^2(V+W)\dd r$, so $V=W$. 
	
	(c) Let $L_V\varphi=\lambda\varphi$ be the principal eigenpair. Since $L_VV=\mu V^2$, self-adjointness gives
	\begin{equation*}
			\lambda\int\varrho V\varphi\dd r=\mu\int\varrho V^2\varphi\dd r>0.
	\end{equation*}
	Thus $\lambda>0$. Smooth dependence follows from the implicit function theorem after mapping an annulus to a fixed interval or a ball to a fixed ball. The  Fredholm theory yields the invertibility of  $L_V$.
\end{proof}

Define positive numbers
\begin{equation}\label{eq-d-l}
	d_1=1,\qquad d_{\ell+1}=d_\ell\frac{\alpha_{\sigma(\ell),\sigma(\ell+1)}}{\alpha_{\sigma(\ell+1),\sigma(\ell)}}.
\end{equation}
\revised{We use these factors to combine the shell solutions into one signed
function. Define $Y$ and $W$ on $[0,R]$ by}  
\begin{equation}\label{eq-Y}
	Y(r): =(-1)^{\ell-1}d_\ell V_\ell(r),\quad r\in I_\ell;
	\qquad W(r): =\frac{|Y|}{d_\ell}=V_\ell(r) \quad r\in I_\ell.
\end{equation}
The relation $d_\ell p_\ell=d_{\ell+1}q_\ell$ makes $Y'$ continuous at every interface point $s_l$ \revised{by \eqref{flux-condi}. Since $Y=0$} there and $Y''=-(n-1)Y'/s_\ell$, $Y''$ is continuous as well. In particular,
\begin{equation}\label{Y-ode}
	-\D Y=\mu Y-\frac\mu{d_\ell}|Y|Y=\mu(1-W)Y\quad\text{on }I_\ell,
	\qquad Y'(0)=0,\quad Y(R)=0.
\end{equation}
The function $W$ is continuous and positive except at the finitely many interfaces and at $R$.

\begin{proposition}\label{morse}
	For $W$ in  \eqref{eq-Y}, define
	\begin{equation*}
			A_0=-\D-\mu+\mu W,\qquad L_W=-\D-\mu+2\mu W
	\end{equation*}
	with regular radial condition at $0$ and zero value at $R$. Then $L_W$ has exactly $m-1$ negative eigenvalues and no zero eigenvalue. More precisely,
	\begin{equation*}
		\lambda_{m-1}(L_W)<0<\lambda_m(L_W)\quad(m\ge2).
	\end{equation*}
	For $m=1$ one has $\lambda_1(L_W)>0$.
\end{proposition}
\begin{proof}
	Since $A_0Y=0$ and $Y$ has exactly $m-1$ interior zeros, Sturm oscillation yields $\lambda_m(A_0)=0$. Strict potential monotonicity, applied to $L_W=A_0+\mu W$, gives
	\begin{equation}\label{positive-m}
		\lambda_m(L_W)>0.
	\end{equation}
	\revised{To prove the other inequality, we use pieces of $Y'$ between
successive critical points as test functions for the quadratic form.} 
	On a positive shell, $(\varrho V_\ell')'=-\mu\varrho V_\ell(1-V_\ell)<0$. Thus $V_1'<0$ for $r>0$ in the first shell, and each later shell has exactly one critical point $t_\ell\in I_\ell$. Consequently $Y'$ has precisely the $m-1$ positive zeros $t_2,\ldots,t_m$. Put $t_1=0$. Differentiating \eqref{Y-ode} in each shell gives
	\begin{equation}\label{Yprime}
		L_W Y'=-\frac{n-1}{r^2}Y'.
	\end{equation}
	For $\ell=1,2, \cdots,m-1$, define $\psi_\ell=Y'$ on $(t_\ell,t_{\ell+1})$ and zero elsewhere. In particular, $Y'(r)=O(r)$ at the center.  
	Multiplying \eqref{Yprime} by $\psi_\ell$ and integration by parts gives
	\begin{equation*}
	\int_0^R\varrho\big[ \psi_\ell'^2+(-\mu+2\mu W)\psi_\ell^2\big] \dd r
		=-(n-1)\int_{t_\ell}^{t_{\ell+1}}r^{n-3}Y'^2\dd r<0.
	\end{equation*}
	\revised{Each piece vanishes at its endpoints and extends by zero to a
function in the radial form domain. Their interiors are disjoint,
so the cross terms in the quadratic form vanish. The span
is therefore $(m-1)$-dimensional and the form is strictly negative
on every nonzero element of the span.} The min--max principle gives $\lambda_{m-1}(L_W)<0$. Together with \eqref{positive-m}, this proves the assertion. For $m=1$, the conclusion holds from Lemma \ref{scalar-radial}.
\end{proof}

\begin{corollary}\label{transmission}
	Suppose $\Phi_\ell$ satisfies $L_{V_\ell}\Phi_\ell=0$ in $I_\ell$, $\Phi_1'(0)=0$, $\Phi_m(R)=0$, and
	\begin{equation}\label{transmission-r}
		\alpha_{\sigma(\ell),\sigma(\ell+1)}\Phi_{\ell+1}(s_\ell)=-\alpha_{\sigma(\ell+1),\sigma(\ell)}\Phi_\ell(s_\ell),\qquad
		\alpha_{\sigma(\ell),\sigma(\ell+1)}\Phi_{\ell+1}'(s_\ell)=-\alpha_{\sigma(\ell+1),\sigma(\ell)}\Phi_\ell'(s_\ell).
	\end{equation}
	Then every $\Phi_\ell$ vanishes.
\end{corollary}
\begin{proof}
	The piecewise function $Z=(-1)^{\ell-1}d_\ell\Phi_\ell$ is $C^1$ by \eqref{eq-d-l} and \eqref{transmission-r}, and satisfies $L_W Z=0$ with the correct endpoint conditions. Proposition~\ref{morse} gives $Z=0$.
\end{proof}

\revised{We now fix the weights $d_\ell$ and choose the value at the
center as a parameter. The initial conditions are}
\begin{equation}\label{shooting-start}
	Y(0;t)=t,\qquad Y'(0;t)=0,\qquad 0<t<1.
\end{equation}
On its first nodal interval, solve $-\D Y=\mu Y-\mu|Y|Y/d_1$. After its first simple zero, use $d_2$; after the next zero, use $d_3$, and so on. Only the coefficient of the quadratic term changes. Values and first derivatives are matched at every crossing. Let $\mathcal T_m\subset(0,1)$ be the set of $t$ such that $m$ finite simple zeros exist, and denote the successive positive zeros by  $b_1(t),b_2(t), \ldots,b_m(t)$. 
\revised{A zero at a positive radius is simple by uniqueness for the
initial-value problem. On a completed nodal interval the function
$|Y|/d_l$ solves the positive logistic equation with the appropriate
endpoint conditions. In particular $0<|Y|/d_l<1$ there by
Lemma~\ref{scalar-radial}. At the center the radial integral formula
gives local existence and smooth dependence on the initial value by a
contraction on a sufficiently short interval, since the integral
operator has norm $O(r^2)$. Away from the center one uses ordinary
differential equation theory. The implicit function theorem at each
simple zero then gives openness of $\mathcal T_m$ and $C^1$
dependence of each $b_j$ and each successive piece on $t$.}

\begin{lemma}\label{return}
	For $t\in\mathcal T_m$, one has $b_m'(t)>0$.
\end{lemma}
\begin{proof}
\revised{Fix $t\in\mathcal T_m$ and restrict $Y(\cdot;t)$ to $[0,b_m(t)]$.
It solves \eqref{Y-ode} there with intervals
$(b_{l-1}(t),b_l(t))$, $b_0(t)=0$, and
$W=|Y(\cdot;t)|/d_l$ on the $l$th interval.}
	Consequently, $Z=\partial_t Y$ solves
	 \begin{equation*}
	 	L_W Z=0,\qquad Z(0)=1,\quad Z'(0)=0,   
	 \end{equation*}
	\revised{and is $C^1$ across the moving zeros. Indeed, differentiating
the matching identities for $Y$ and $Y'$ along a zero cancels the
terms containing the derivative of the zero because both $Y'$
and $Y''$ agree there. By Proposition} \ref{morse},   $L_W$ has exactly $m-1$ negative eigenvalues and no zero eigenvalue, and then  \revised{the Sturm oscillation theorem for the regular solution at spectral
parameter zero gives exactly $m-1$ zeros of $Z$ before $b_m(t)$}, and the sign of $Z(b_m(t))$ is $(-1)^{m-1}$.
	Moreover, notice that $\sgn Y'(b_m(t))=(-1)^m$. Differentiating $Y(b_m(t);t)=0$  yields
	\begin{equation*}
		b_m'(t)=-\frac{Z(b_m(t);t)}{Y'(b_m(t);t)}>0.
	\end{equation*}
%	The same argument applies to every earlier completed return.
\end{proof}

\begin{lemma}\label{range}
	There exists $t_m^*\in(0,1]$ such that
	\begin{equation*}
			\mathcal T_m=(0,t_m^*),\qquad
		\lim_{t\downarrow0}b_m(t)=\frac{j_{\upnu,m}}{\sqrt\mu},\qquad
		\lim_{t\uparrow t_m^*}b_m(t)=+\infty.
	\end{equation*}
\end{lemma}
\begin{proof}
	\revised{Write $Y=t\varphi_t$. In the equation for $\varphi_t$, the quadratic
term has coefficient $\mu t/d_l$. Starting from
$\varphi_t(0)=1$, continuous dependence gives convergence on bounded
intervals of the first branch to the regular solution $\varphi$ of}
	\begin{equation*}
			-\D\varphi=\mu\varphi,\qquad\varphi(0)=1,\quad\varphi'(0)=0.
	\end{equation*}
	\revised{Thus $\varphi(r)=C_\mu r^{-\upnu}J_{\upnu}(\sqrt{\mu}\,r)$ with
$C_\mu=2^{\upnu}\Gamma(\upnu+1)\mu^{-\upnu/2}$, so that
$\varphi(0)=1$. Its positive zeros are simple and equal to
$j_{\upnu,j}/\sqrt\mu$.} \revised{Because the first zero of $\varphi$ is simple, the first
zero of $\varphi_t$ exists for small $t$ and converges to it.
After switching to $d_2$, the coefficient of the quadratic term
still tends to zero, and the crossing data converge. Repeating
this argument up to the $m$th zero proves both that
$(0,\delta)\subset\mathcal T_m$ for some $\delta>0$ and that
$b_m(t)\to j_{\upnu,m}/\sqrt\mu$.}
	
	Suppose $t$ stays in a compact subset of $(0,1)$ and $b_m(t)\le B$. Then
	\begin{equation*}
			-\D Y=q_t(r)Y,\qquad 0\le q_t(r)=\mu(1-|Y|/d_\ell)\le\mu.
	\end{equation*}
	On a sufficiently small fixed interval $[0, \delta_0]$, $Y(r)=t+O(r^2)$ is uniformly positive. 
	On $[\delta_0,B]$, one has
	\begin{equation*}
			\frac{\mathrm{d}}{\mathrm{d}r} [Y^2 +(Y')^2]
=
2Y Y' + 2Y'\big[-q_t(r) Y - \frac{n-1}{r} Y' \big] 
\ge 
- \big[
 \revised{1+\mu+\frac{2(n-1)}{\delta_0}}
 \big]  [Y^2 +(Y')^2].
	\end{equation*}
	\revised{The initial values range in a compact subset of $(0,1)$, so
$Y(\delta_0)$ is bounded below by a positive constant. Gronwall's
inequality now bounds $|Y|+|Y'|$ from below on
$[\delta_0,b_m(t)]$. In particular, the derivatives at zero crossings
are bounded away from zero. The same differential equation bounds
$Y'$ and $Y''$ uniformly on that interval, by the uniform bound on
$Y$ and integration from $\delta_0$. At consecutive zeros the
derivatives have opposite signs, so the uniform bound on $Y''$
keeps those zeros apart. The first zero stays away from $0$ by
the initial positivity. Along a convergent sequence of parameters,
continuous dependence up to each simple crossing then gives
successively $m$ distinct simple zeros for the limiting parameter.}
	
	 Let $(a,b)$ be any component of $\mathcal T_m$. By Lemma~\ref{return}, $b_m(t)$ is increasing on it. If $a>0$, then $b_m(t)$ is bounded above as $t\downarrow a$ by its value at any fixed point of the component. The preceding compactness argument extends the returns through $a$, contradicting the choice of the component. Thus every component has left endpoint $0$, and then $\mathcal T_m=(0,t_m^*)$.
 The monotone function $b_m$ has a limit there. If the limit were finite and $t_m^*<1$, compactness would extend the domain beyond $t_m^*$. If $t_m^*=1$, continuous dependence in the first shell gives $Y(\cdot;t)\to1$ on every fixed finite interval, because the equation has the regular solution $Y\equiv1$ at initial value $1$. Its first zero cannot remain bounded. In either case, a finite limiting $m$-th zero is impossible.
\end{proof}

\subsection{Maximum principles for linear analysis and positivity}
\begin{lemma}\label{global-barrier}
	Let $|q_\eps|\le Q$ uniformly bounded, and $S_\eps\ge\xi\min\{\eps,R-r\}$ on $[0,R]$, with $\xi>0$ fixed. For small $\eps$, the operator
	\begin{equation*}
			B_\eps=-\D+\eps^{-3}S_\eps-q_\eps,
		\qquad v'(0)=0,\quad v(R)=0,
	\end{equation*}
	has the maximum principle and is invertible. Moreover,
	\begin{equation}\label{scalar-inv}
		\|B_\eps^{-1}f\|_\infty\le C\eps^2\|f\|_\infty,\qquad
		\|B_\eps^{-1}(\eps^{-3}S_\eps g)\|_\infty\le C\|g\|_\infty.
	\end{equation}
\end{lemma}
\begin{proof}
	Set
	\begin{equation*}
		b_\eps(r)=1-\frac{\cosh(r/\eps)}{\cosh(R/\eps)}.
	\end{equation*}
	It has $b_\eps'(0)=0$, $b_\eps(R)=0$, $0\le b_\eps\le1$, and
	\begin{equation*}
			-\D b_\eps=\frac{\eps^{-2}\cosh(r/\eps)+(n-1)(\eps r)^{-1}\sinh(r/\eps)}{\cosh(R/\eps)}
		\ge
		\eps^{-2}\frac{\cosh(r/\eps)}{\cosh(R/\eps)}.
	\end{equation*}
	The expression at $r=0$ is interpreted by its regular limit,
	and
	\begin{equation*}
			-\D b_\eps(0) =
			\lim\limits_{r\rightarrow 0^+} \frac{\eps^{-2}\cosh(r/\eps)+(n-1)(\eps r)^{-1}\sinh(r/\eps)}{\cosh(R/\eps)} 
			= \frac{n \eps^{-2}}{\cosh(R/\eps)}.
	\end{equation*}
	\revised{For $R-r\le\eps$}, the ratio $\frac{\cosh(r/\eps)}{\cosh(R/\eps)} \ge \frac{1}{2\mathrm{e}}$. For $R-r\ge\eps$, $b_\eps\ge1-2e^{-1}>0$ and $\eps^{-1}S_\eps\ge\xi$. Thus
	\begin{equation*}
			(-\D+\eps^{-3}S_\eps)b_\eps\ge c_0\eps^{-2},\qquad
		B_\eps b_\eps\ge\tfrac {c_0}{2} \eps^{-2}, \quad c_0=\min \{ \frac{1}{2\mathrm{e}}, \xi(1-\frac{2}{\mathrm{e}}) \}
	\end{equation*}
	for small $\eps$.
	Let $h_\eps:=b_\eps+1$, then  $h'_\eps (0)=0$, $h_\eps (R)=1$, $1\le h_\eps \le 2$, and $B_\eps h_\eps \ge \tfrac {c_0}{2}\eps^{-2}-2Q >0 $ for small $\eps$. 
	
	Now suppose $B_\eps v\ge 0$, $v'(0)=0$ and $v(R)\ge0$, but $\min v<0$. \revised{Set $w=v/h_\eps$. Since $h_\eps\ge1$, we have} $w\in C^2([0,R])$. Moreover,
	\begin{equation*}
		w'(0)=\frac{v'(0)h_\eps(0)-v(0)h_\eps'(0)}{h_\eps(0)^2}=0,
		\qquad
		\revised{w(R)=v(R)\ge0}.
	\end{equation*}
	Furthermore,  we see that 
	\begin{equation}\label{w-eq-quotient}
	0\le 	 \frac{B_\eps v}{h_\eps}
		=-w''-\left(\frac{n-1}{r}
		+2\frac{h_\eps'}{h_\eps}\right)w'
		+a_\eps w, \quad a_\eps(r)=\frac{B_\eps h_\eps(r)}{h_\eps(r)}>0.
	\end{equation}
	Assume that $w(r_0)=\min w<0$, there are two cases:   (a) if $r_0\in (0,R)$, at the minimum $r_0$, one has $w'(r_0)=0$, and $w''(r_0)\ge 0$, but then 
	\begin{equation*}
	0\le 	 \frac{B_\eps v(r_0)}{h_\eps(r_0)}
		=-w''(r_0)-\left(\frac{n-1}{r_0}
		+2\frac{h_\eps'(r_0)}{h_\eps(r_0)}\right)w'(r_0)
		+a_\eps(r_0) w(r_0)  =-w''(r_0)+a_\eps(r_0) w(r_0)<0,
	\end{equation*}
	which is a contradiction.
	\revised{(b) If $r_0=0$, then since} $w'(0)=0$, it holds that $\frac{n-1}{r} w'(r)\rightarrow (n-1)w''(0)$ as $r\rightarrow 0^+$. And the minimum property yields
	\begin{equation*}
		0\le w(r)-w(0) = \frac{w''(0)}{2} r^2 +o(r^2). 
	\end{equation*}
	Consequently, $w''(0)\ge 0$. Hence by taking the limit $r\rightarrow 0^+$ in \eqref{w-eq-quotient}, we obtain
	\begin{equation*}
		0\le -n w''(0) +a_\eps (0) w(0) <0, 
	\end{equation*}
\revised{This contradiction also excludes a negative minimum at the center.
Hence the maximum principle holds on the whole radial interval.} 	
 Kernel triviality and Fredholm theory give invertibility.
	
	Comparison with $C\eps^2\|f\|_\infty b_\eps$ gives the first estimate in \eqref{scalar-inv}. For the second, choose $C_0$ so that
	\begin{equation*}
			B_\eps(1+C_0\eps^2b_\eps)\ge\eps^{-3}S_\eps.
	\end{equation*}
	This follows from $B_\eps1=\eps^{-3}S_\eps-q_\eps$. The barrier is positive at $R$ and has zero derivative at $0$, so comparison gives the second estimate.
\end{proof}

\begin{lemma}\label{gap-left-out-barrier}
	\revised{Assume $|q_\eps|\le Q$ uniformly. There is $\zeta_2>0$ such that
the following estimates hold for the operator
$-\D+\eps^{-3}S_\eps-q_\eps$ for every $0\le\zeta\le\zeta_2$
and all sufficiently small $\eps$.}
	
	(a)  On $[a+\eps,R]$, where $a>0$ is fixed and $S_\eps\ge\xi\min\{\eps,R-r\}$, use
	\begin{equation*}
		(-\D+\eps^{-3}S_\eps-q_\eps)B^{\mathrm{out}}
	\ge c\eps^{-2}e^{-\zeta(r-a)/\eps}, \quad \mbox{where }\ 	B^{\mathrm{out}}(r):=e^{-\zeta(r-a)/\eps}(1-e^{-(R-r)/\eps}).
	\end{equation*}
	
	(b) On $[a+\eps,b-\eps]$, with fixed $0<a<b<R$ and $S_\eps\ge\xi\eps$, use
	\begin{equation*}
		(-\D+\eps^{-3}S_\eps-q_\eps)B^{\mathrm{gap}}
	\ge c\eps^{-2}B^{\mathrm{gap}},  \quad \mbox{where }\ 	B^{\mathrm{gap}}(r):=e^{-\zeta(r-a)/\eps}+e^{-\zeta(b-r)/\eps}.
	\end{equation*}
	
(c) On $[0,a-\eps]$, with $a>0$ fixed and $S_\eps\ge\xi\eps$, use
	\begin{equation*}
	(-\D+\eps^{-3}S_\eps-q_\eps)B^{\mathrm{cen}}
	\ge 
	c\eps^{-2}B^{\mathrm{cen}}, \quad \mbox{where }\ 
		B^{\mathrm{cen}}(r):=e^{-\zeta(a-r)/\eps}+e^{-\zeta(a+r)/\eps}
		=2e^{-\zeta a/\eps}\cosh(\zeta r/\eps).
	\end{equation*}
	It has zero derivative at $0$, and $	e^{\zeta(a-r)/\eps}B^{\mathrm{cen}}(r)\le2$. 
	All the corresponding scalar boundary problems have the maximum principle.
\end{lemma}
\begin{proof}
	On regions away from zero, $(n-1)/r$ is bounded. For the outer barrier, put $t=(R-r)/\eps$. Differentiation gives
	\begin{equation*}
		\eps^2e^{\zeta(r-a)/\eps}(-\D+\eps^{-3}S_\eps-q_\eps)B^{\mathrm{out}} 
	 \ge(1-2\zeta-C\eps)e^{-t}
		+[\xi\min\{1,t\}-\zeta^2-C\eps\zeta-Q\eps^2](1-e^{-t}).
	\end{equation*}
	At $\zeta=0$, the two positive terms have a strictly positive infimum for $t\ge0$. Choose $\zeta$ small and then $\eps$ small. This proves the outer estimate.
	
	 For the gap region, direct calculations yield
	 \begin{equation*}
	 	(-\D+\eps^{-3}S_\eps-q_\eps)B^{\mathrm{gap}} 
	 	\ge 
	 	\eps^{-2} B^{\mathrm{gap}} \big[
	 	-\zeta^2 -\frac{n-1}{r} \zeta \eps +\xi -\eps^2 Q
	 	\big],
	 \end{equation*}
	 giving the assertion.
	
	At the center region, using $\tanh t\le t$ for $t\ge0$ gives
	\begin{equation*}
		\D B^{\mathrm{cen}} = B^{\mathrm{cen}} \big[
		\revised{\zeta^2\eps^{-2}+\frac{n-1}{r}\frac{\zeta}{\eps}\tanh\!\left(\frac{\zeta r}{\eps}\right)}
		\big]
		\le 
		B^{\mathrm{cen}} \big[
		\revised{\zeta^2\eps^{-2}+(n-1)\zeta^2\eps^{-2}}
		\big]
		\le n\zeta^2\eps^{-2}B^{\mathrm{cen}}.
	\end{equation*}
	Choose also $n\zeta^2<\xi/2$. This proves the central estimate. \revised{For the gap and center problems the potential is positive for
small $\varepsilon$, so the maximum principle follows at a negative
minimum, using the regular radial limit at $0$ in the center problem.
For the outer problem use $B^{\mathrm{out}}$ at $\zeta=0$ plus a
small positive constant. This is strictly positive on the closed
interval and remains a strict supersolution, as in the proof of
Lemma~\ref{scalar-barriers-1}.}
\end{proof}

\section*{Acknowledgments}
The authors used ChatGPT to assist in simplifying and refining the presentation of the proofs. 
Z. Li is supported by the National Natural Science Foundation of China (NSFC)
[Grant No. 12501138], and funded by Science Research Project of Hebei Education Department [Grant No.
BJ2026100],  the Hebei Province Yanzhao Golden Stage Talent Gathering Program Key Talent Project
(Study Abroad Returning Platform) [Grant No. B2024018], the China Scholarship Council (CSC) [Grant No. 202508130126], and the 111 Center [Grant No. D26018].

%\bibliographystyle{plain}
%\bibliography{RefData}

%\end{document}

%\iffalse

%\fi
\end{document}